# An exact algorithm for the pickup and delivery problem with time windows and requests split across multiple stacks

Ali Mehsin Alyasiry, Michael Forbes*, Ella Wang

*School of Mathematics and Physics, Faculty of Science, The University of Queensland, Brisbane, Australia*

**Abstract**

This paper addresses an extension of the pickup and delivery problem with time windows and multiple stacks (PDPTWMS). We call the new extension the pickup and delivery problem with time windows and requests split across multiple stacks (PDPTWRSMS). Applications of the PDPTWRSMS arise in less-than-truckload transportation systems where loads are packed as pallets or other standardised containers. In the PDPTWRSMS, the vehicle's loading space is divided into compartments of finite capacity. Loading and unloading in each compartment must obey the last-in-first-out (LIFO) policy. Moreover, a request can be split among all compartments, and a single request can also exceed a single compartment's capacity. Thus, each request is only limited to the vehicle's capacity. Computational experiments show that applying this extension can lead to substantial routing cost reductions (the total distance travelled and/or the number of vehicles used), and, compared with versions of the PDPTWMS where compartments are not restricted to a LIFO policy, splitting is worth more once the vehicle has three or more compartments, though less when it has two deep ones. We apply a modified version of the recently proposed method of fragments. Results confirm that this approach can solve many PDPTWRSMS instances and significantly outperforms the current state-of-the-art method for the PDPTWMS.

## 1. Introduction

In a one-to-one pickup and delivery problem with time windows (PDPTW) vehicles are routed to serve given requests each of which consists of a pickup and a corresponding delivery. For each request, the pickup must precede the delivery (precedence) and both must be performed by the same vehicle (pairing). Routes must respect precedence, pairing, vehicle-capacity, and time-window constraints. A set of feasible routes is a solution to the PDPTW if each request is performed by exactly one of the routes.

A further specialisation of the PDPTW is the pickup and delivery problem with time windows and last-in-first-out (LIFO) loading (PDPTWL). In this problem loading and unloading an item must obey the LIFO loading policy, meaning loaded items are treated as if they are placed on a stack in the vehicle, see (Cherkesly et al., 2015). This will avoid extra cost associated with rearranging items inside the vehicle's loading compartment. LIFO loading is very common in the transportation of pallets, animals, heavyweight goods and hazardous materials. The problem is also seen in *Automated Guided Vehicle* systems that are often used in sea ports and warehouses.

The pickup and delivery problem with time windows and multiple stacks (PDPTWMS), introduced in Cherkesly et al. (2016) is a generalisation of the PDPTWL for vehicles with multiple load compartments. In the PDPTWMS, loading and unloading an item in each compartment must obey the LIFO loading policy. This is called the multi-stack policy. A vehicle route is feasible for the PDPTWMS if it is feasible for the PDPTW and respects the multi-stack policy. The objective of the PDPTWMS is to satisfy all requests with a least-cost routing. Each route incurs a travel cost and a fixed cost (vehicle cost). We assume an unlimited supply of identical vehicles. Each vehicle consists of $k = |S|$ independently accessible identical compartments

*Corresponding author

*Email addresses:* a.alyasiry@uq.net.au (Ali Mehsin Alyasiry), m.forbes@uq.edu.au (Michael Forbes), ella.wang@uq.edu.au (Ella Wang)

of finite capacity $\hat{Q}$, where $S = \{1, ..., k\}$ is the set of the compartments in any vehicle. The total vehicle capacity is $Q = k \times \hat{Q}$.

The only works in the literature studying the PDPTWMS are that by Cherkesly et al. (2016) and Cherkesly and Gschwind (2022). Cherkesly et al. (2016) present two branch-price-and-cut algorithms. The first algorithm solves the pricing problem by fully enforcing the multi-stack policy. The second algorithm relaxes the definition of a feasible route in the pricing problem and infeasible path inequalities are added to the relaxed master problem if infeasible routes are used in a linear relaxation solution. Instances with up to 75 requests are solved to optimality within two hours of computation time.

Cherkesly and Gschwind (2022) extend the PDPTWMS to the pickup and delivery problem with time windows, multiple stacks and handling operations (PDPTWMS-H). This extension still requires that a request be kept together on one stack, but allows handling operations at delivery locations to rearrange the stacks.

A related line of work divides the vehicle's loading space into compartments but does not require the compartments to be operated under a LIFO policy. Aerts-Veenstra et al. (2024) introduce the pickup and delivery problem with time windows and multiple compartments (PDPTWMC) and study three compartment-related attributes: compartment capacity flexibility, where the capacity of each compartment is either fixed or chosen before the route subject to a minimum and a maximum; item-to-compartment flexibility, which specifies which requests may be loaded in which compartments; and item-to-item compatibility, where incompatible requests cannot be in the same compartment at the same time. They propose a branch-price-and-cut algorithm in which the pricing problem is solved by a unified bidirectional labelling algorithm, and they introduce two symmetry-reduction techniques, one which suppresses extensions to interchangeable empty compartments and one which allows dominance to hold under a permutation of comparable compartments. Instances with up to 75 requests are solved to optimality within one hour of computation time.

Aerts-Veenstra et al. (2025) extend this family to the pickup and delivery problem with time windows, multiple compartments, incompatibility constraints and cleaning operations (PDPTWMCI). They distinguish four types of incompatibility between a pair of requests, according to whether the conflict applies to a single compartment or to the whole vehicle, and whether it applies only while both requests are on board or also to a request carried subsequently. Incompatibilities of the latter kind can be removed by cleaning a subset of empty compartments at a dedicated cleaning station, at a cost of a fixed time plus a time proportional to the number of compartments cleaned. The pricing problem is again solved by labelling, with each label carrying, for each compartment, the set of requests that the compartment can no longer accept.

The loading policy of these problems is different from the multi-stack policy considered here. A compartment imposes no order on the requests it holds, so the restrictions on a route arise from the compatibility of a request with a compartment or with the requests already on board rather than from LIFO loading, and the cleaning operations of the PDPTWMCI play a role analogous to that of the handling operations of Cherkesly and Gschwind (2022), in that both make feasible a route that would otherwise not be, at the cost of additional time spent at an intermediate vertex. Compartment capacity flexibility does allow a compartment to be as large as the vehicle, so that the size of a request need not be bounded by an equal division of the vehicle capacity. The capacities are, however, chosen before the route and remain fixed throughout it, and a request is always assigned to a single compartment. As Aerts-Veenstra et al. (2024) observe, the compartment-related attributes they study have not been considered in the multiple stacks literature.

In the literature several papers considered the single-vehicle pickup and delivery problem with LIFO loading or first-in-first-out (FIFO) loading, see, e.g., (Carrabs et al., 2007; Cordeau et al., 2010a, 2010b; Erdoğan et al., 2009, Li et al., 2011), and the pickup and delivery problem with time windows and LIFO loading (PDPTWL), see, e.g., (Alyasiry et al., 2019, Cherkesly et al., 2015), for the single-stack case.

For the multiple stacks case there are papers on the single-vehicle pickup and delivery problem with multiple stacks, see, e.g., (Côté et al., 2012a, 2012b; Pereira and Urrutia, 2018, Sampaio and Urrutia, 2017), the double travelling salesman problem with multiple stacks (DTSPMS), see, e.g., (Alba Martínez et al., 2013, Carrabs et al., 2013, Felipe et al., 2009, Lusby et al., 2010, Petersen et al., 2010, Petersen and Madsen, 2009), and the double vehicle routing problem with multiple stacks (DVRPMS) (Iori and Riera-Ledesma, 2015).

To the best of our knowledge, all previous methods in the literature that have dealt with multiple stacks routing have assumed that a request cannot be split among different stacks and that the size of a request cannot exceed the stack capacity. Requests are not split among compartments in the multi-compartment literature either, although Aerts-Veenstra et al. (2024) allow the capacity of each compartment to be chosen, within given bounds, before the route begins.

In this paper we introduce an extension to the PDPTWMS. We allow the size of a request to exceed the stack capacity and only be limited by the vehicle capacity, and we allow a request to be split among all stacks. In the vehicle routing problem literature splitting often refers to splitting a request to transport it on two or more vehicles, for more details see (Archetti and Speranza, 2012). In contrast, the splitting we consider is among the compartments of one vehicle. Therefore, we name the new variant the pickup and delivery problem with time windows and requests split across multiple stacks (PDPTWRSMS). Computational experiments show a reduction in the total travelled distance and/or a reduction in the number of vehicles used in many instances when request splitting is allowed.

### *1.1. Contributions*

In this paper, we extend the work of Alyasiry et al. (2019) for the PDPTWL and adapt that method to the PDPTWMS and the PDPTWRSMS. Alyasiry et al. (2019) introduced a branch-and-cut approach based on fragments of requests. A fragment is a sequence of correctly paired and scheduled pickup and delivery vertices starting and ending with an empty vehicle. Using fragments, they formulate a relaxed network flow model with side constraints.

We provide a general purpose framework to generate all feasible fragments for a range of PDPTW problems, including the PDPTW, PDPTWL, PDPTWMS, PDPTWRSMS and the PDPTWMS-H with zero handling time. In all cases we generate fragments while ensuring that we can assign pickups to compartments so as to respect the relevant multi-stack policy. Such an assignment of pickups is called a loading plan.

Following Rist and Forbes (2021) and Rist and Forbes (2022), we extend the fragments to include the next pickup or the depot. We delete dominated extended fragments and formulate an IP to chain the fragments to form routes using a branch-and-cut approach similar to that of Rist and Forbes (2022).

We introduce an integer programming (IP) formulation to assign requests to stacks to produce a feasible multi-stack loading plan. In addition, we use a number of techniques to minimise the number of stack assignment IPs solved and show how to produce load plans using dynamic programming when requests cannot be split across stacks.

Further we conduct extensive computational tests to evaluate the benefits of request splitting. In particular we compare the PDPTWMS, the PDPTWRSMS and the PDPTWMS-H under handling policy H1 with zero handling time. Computational results show that our method outperforms the current state-of-the-art method for the PDPTWMS and solves to optimality many instances that were previously unsolved.

### *1.2. Outline*

The remainder of this paper is organised as follows. The next section provides notation. Section 3 formally defines the PDPTWMS and PDPTWRSMS and describes the difference between them. In Section 4 we define the fragments and we present procedures for generating loading plans (assigning pickups to compartments) including an IP formulation. In Section 5 we present the mathematical formulation used for all PDPTW variants and describe solution acceleration techniques. In Section 6 we report computational results. In Section 7 we examine the optimal solutions themselves, and compare the saving from request splitting to the saving from relaxing the LIFO constraint on stacks. Conclusions are presented in Section 8.

## 2. Notation

All PDPTW variants considered can be modelled using a directed graph (network) $G = (V, A)$, where $V$ represents the set of nodes or vertices and $A$ represents the set of arcs. The set of vertices for a problem with $n$ requests can be characterised as $V = \{0\} \cup P \cup D$, where vertex 0 represents the depot; $P = \{1, \cdots, n\}$ is a set of pickup vertices; and $D = \{n + 1, \cdots, 2n\}$ is a set of delivery vertices. The set of arcs $A$ can be characterised as $A = \{(i, j) : i, j \in V, i \neq j\}$. Infeasible arcs are eliminated using the rules presented in Dumas et al. (1991). Each pickup vertex $i \in P$ has a corresponding delivery vertex $n + i \in D$. Due to the pairing conditions between the pickup and the delivery vertices, the set of pickup vertices $P$ can also represent the set of requests.

For each vertex $i \in V$ there is a load (the size of the request) $q_i > 0$ if $i \in P$ and an unload $q_i = -q_{i-n}$ if $i \in D$, and a time window $[e_i, l_i]$, where $e_i$ and $l_i$ represent the earliest and the latest time the service at vertex $i$ can begin. Time windows are tightened using the rules presented in Dumas et al. (1991). The load at the depot $q_0 = 0$ and the time window at the depot is $[e_0 = 0, l_0 = L]$, where $L$ is the length of the planning horizon.

The travel cost from vertex $i$ to vertex $j$ is denoted by $c_{ij}$, and the travel time is denoted by $t_{ij}$ and includes a non-negative service time duration at vertex $i$. We assume that travel costs and travel times satisfy the triangle inequality.

A series of $m$ vertices $\langle v_1, v_2, ..., v_m \rangle$ is a fragment if $\langle 0, v_1, v_2, ..., v_m, 0 \rangle$ is a valid route, satisfying time window, capacity, pairing and precedence constraints, and moreover the vehicle is never empty. That is $\sum_{i=1}^{m'} q_{v_i} > 0$ for all $1 < m' < m$. Any legal route can be formed by chaining one or more fragments together and adding the depot to each end.

**Definition 2.1.** An extended fragment is a fragment with an additional pickup vertex, or the depot, added to the end of the fragment. A fragment can always be extended to the depot. A fragment $f = \langle v_1, v_2, ..., v_m \rangle$ can be extended to a pickup vertex $p \in P$ if $p \notin f$ and $\langle 0, v_1, v_2, ..., v_m, p, p + n, 0 \rangle$ is a legal route.

Any legal route can be formed by chaining extended fragments together where the last vertex of one fragment is the same as the first vertex of the next fragment. Only one copy is kept of the repeated vertex.

## 3. PDPTWRSMS vs PDPTWMS vs PDPTWMS-H

Using multiple LIFO stacks increases the loading options compared with a single LIFO stack and may lead to a significant saving, since many loading patterns that are infeasible with one stack become feasible with more. Suppose, for example, that the problem consists of two requests with pickup vertices $p_1$, $p_2$ and delivery vertices $d_1$, $d_2$. The loading and unloading pattern $\langle p_1, p_2, d_1, d_2 \rangle$ is feasible for two or more stacks of sufficient size and infeasible for a single stack of any size, since $d_1$ is reached while request 2 still lies on top of request 1. When the pickups and deliveries of two requests interact in this way, whatever else the route contains, the requests are said to be stack incompatible.

In the PDPTWMS for any route to be feasible it must be feasible for the PDPTW and it must satisfy the capacity and LIFO constraints for each stack (Cherkesly et al., 2016). Also, a request must be placed completely on one stack. However, in real-world situations a request often consists of pallet loads which may be split among all stacks, and the size of any request may exceed the stack capacity and only be limited to the vehicle's capacity.

To illustrate, consider a vehicle with two compartments, each of capacity three, and requests of size two. Any path that begins with three pickups is infeasible for the PDPTWMS, since a compartment of three cannot hold two requests and the two compartments therefore hold two requests between them, yet many such paths are feasible for the PDPTWRSMS. Figure 1 shows one of them, the path $\langle p_1, p_2, p_3, d_3, d_1, d_2 \rangle$, together with the loading plan that makes it feasible. A request of size greater than three would have to be rejected altogether in the PDPTWMS, since no single compartment could hold it.

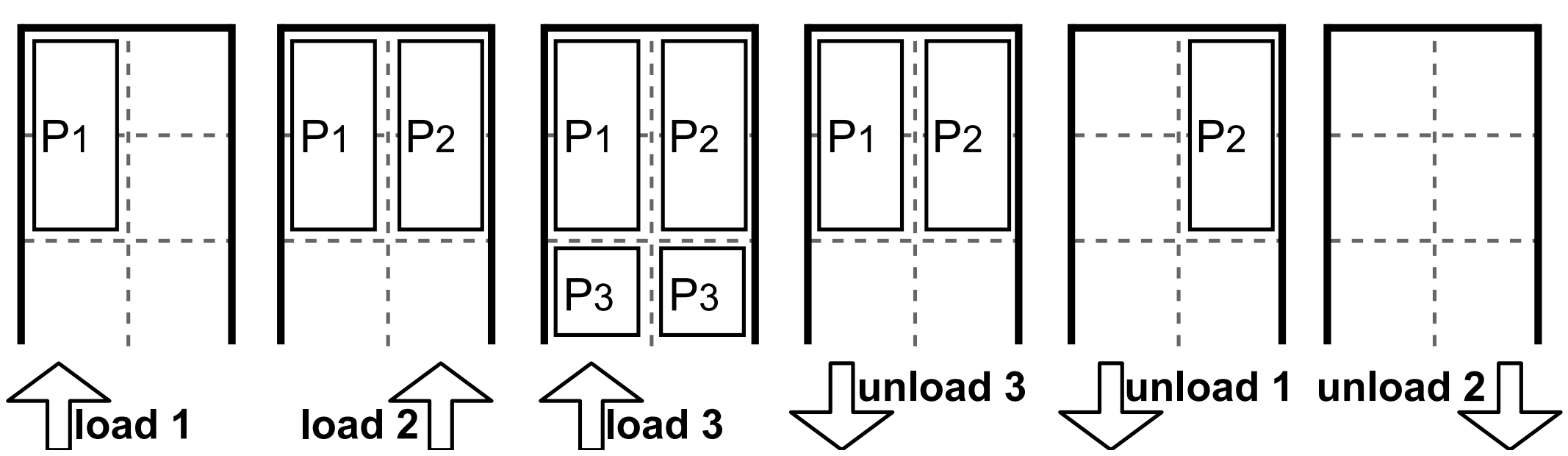

Figure 1: A feasible loading plan for the PDPTWRSMS that is infeasible for the PDPTWMS

More compartments are not always better, even when requests may be split. For example, suppose the vehicle capacity is 12 and the size of each request is four. Ignoring the time windows, the path $\langle p_1, p_2, p_3, d_1, d_2, d_3 \rangle$ is feasible for the three-stack case (the capacity of each compartment is four). Whereas, it is infeasible for the four-stack case even if requests can be split across multiple stacks (the capacity of each compartment is three).

The same can occur between two stacks and three. Suppose the vehicle capacity is 24 and a route consists of four requests of size six, with the loading and unloading pattern $\langle p_1, p_2, p_3, p_4, d_2, d_1, d_4, d_3 \rangle$. With two stacks each compartment holds twelve and the pattern is feasible. Requests 1 and 2 are loaded on one compartment and 3 and 4 on the other, and each pair is unloaded in the reverse of the order in which it was loaded. With three stacks each compartment holds only eight, which is not enough for two requests of size six, so the four requests on board after $p_4$ cannot be placed. Splitting does not repair it as $p_3$ and $p_4$ cannot share stacks with $p_1$ and $p_2$.

### *3.1. PDPTWMS-H*

Cherkesly and Gschwind (2022) introduce the PDPTWMS-H. In this problem requests can be moved within or between stacks at any delivery vertex, with a handling time proportional to the total size of the requests moved. Unless otherwise stated, we only consider the case with zero handling time in this paper. The specific rules around moving requests are defined by the handling policy. Here we consider the policy which corresponds to handling policy H1 in Cherkesly and Gschwind (2022): it allows any compulsory unloading so as to be able to unload the request for the current delivery vertex, with that request's stack reloaded as required. With zero handling time, H1 corresponds to treating the stacks as separate compartments, without imposing a LIFO policy.

In terms of the sets of fragments produced we can define the following relationships: PDPTWMS $\subseteq$ PDPTWRSMS $\subseteq$ PDPTW; and PDPTWMS $\subseteq$ PDPTWMS-H1 $\subseteq$ PDPTW. It is straightforward to produce fragments which have feasible load plans for the PDPTWRSMS but not for the PDPTWMS-H1, and vice versa. Splitting divides a request between stacks but keeps each stack in LIFO order, while H1 keeps every request whole but imposes no order at all. Neither relaxation contains the other, and Section 7.4 compares the savings from each against the PDPTWMS.

## 4. Fragments

Fragments must satisfy pairing, precedence, time window, and vehicle capacity constraints. In addition, fragments must satisfy the multi-stack policy. That is, there must exist an assignment of requests to stacks that respects the capacity of each stack and the LIFO loading policy for each stack. Multi-stack fragments, especially with requests split across stacks, are a non-trivial extension of the PDPTWL fragments introduced by Alyasiry et al. (2019).

For a fragment to be feasible for the PDPTWMS or for the PDPTWRSMS it must first be feasible for the PDPTW and there must exist a loading plan for the fragment that assigns requests to stacks so as to respect both stack capacity and the LIFO policy. More formally, a feasible loading plan for a fragment specifies the quantity of each request assigned to each stack such that:

- Every request is completely assigned to the stacks.
- No pair of stack incompatible requests are assigned to the same stack.
- The total quantity in any stack at any time is less than the stack capacity.

For the PDPTWMS, each request is assigned to exactly one stack.

We denote a fragment by $f$ and the set of all fragments by $\Phi$. We denote by $\theta_f$ the ordered set of requests in $f$, ordered by their appearance in $f$. For all $p \in \theta_f$ we denote by $X_f(p)$ the pickup position of request $p$ in $f$, and $\Gamma_f(p)$ the delivery position of request $p$ in $f$. For example, if $f = \langle p_1, p_2, d_1, p_3, d_3, d_2 \rangle$, then $\theta_f = \{p_1, p_2, p_3\}$, $X_f(p_1) = 1$, $X_f(p_2) = 2$, $X_f(p_3) = 4$, $\Gamma_f(p_1) = 3$, $\Gamma_f(p_2) = 5$, and $\Gamma_f(p_3) = 6$. Here $p_i \in P$ and $d_i = p_i + n$.

Using $X_f(p)$ and $\Gamma_f(p)$ we can compute the set of stack incompatible pairs of requests for fragment $f$. These requests cannot be placed on the same stack because of the LIFO loading policy.

**Definition 4.1.** We say that the requests $p' \in \theta_f$ and $\tilde{p} \in \theta_f$ are *stack incompatible* in fragment $f$ if either of the following conditions are met:

- $X_f(p') < X_f(\tilde{p})$, $\Gamma_f(p') < \Gamma_f(\tilde{p})$ and $X_f(\tilde{p}) < \Gamma_f(p')$;
- $X_f(\tilde{p}) < X_f(p')$, $\Gamma_f(\tilde{p}) < \Gamma_f(p')$ and $X_f(p') < \Gamma_f(\tilde{p})$.

For all $p \in \theta_f$ we denote the set of incompatible requests in $f$ by $I_f(p)$.

### *4.1. Generating Fragments*

We generate all fragments using an efficient recursive algorithm as shown in pseudo-code in Algorithm 1. The procedure *Extend* produces a sequence of vertices to generate a candidate fragment $f$, then the procedure *CanLoad* checks whether the generated fragment $f$ can be feasibly assigned to the vehicle compartments.

*Extend* is a recursive function which takes a partial fragment, the current time and the current load. It first determines if the partial fragment can be completed (i.e. all requests delivered) and, if it can be completed, then tries all fragment extensions to pickups. The construction process dynamically enforces all PDPTW rules – pairing, precedence, time windows and vehicle capacity. Immediately before a candidate fragment is added to the set of feasible fragments we check for a feasible loading plan. Note that we never add pickups to a partial fragment unless the fragment can be feasibly completed, where feasibility includes checking for a loading plan.

**Algorithm 1** Generate all PDPTWMS feasible fragments

**function** EXTEND($f, t', q'$)
    **if** $q' = 0$ **then**
        $S \leftarrow (\{\}, ..., \{\})$ ▷ Empty all stacks
        **if** CANLOAD($f', S$) **then**
            $\Phi \leftarrow \Phi \cup \{f\}$
            **return** `true`
        **else**
            **return** `false`
    `CanComplete` $\leftarrow$ `false`
    $r_m \leftarrow$ the last vertex in $f$
    **for** $r \in P \cap f$ **do** ▷ Try to completely deliver everything first
        **if** $(r + n) \notin f$ **then**
            $t'' \leftarrow \min\left(t' + t_{r_m,(r+n)}, e_{r+n}\right)$ ▷ Arrival time at $(r + n)$
            **if** $t'' \leq l_{r+n}$ and EXTEND($f + \langle r + n \rangle, t'', q' - q_r$) **then**
                `CanComplete` $\leftarrow$ `true`
    **if** not `CanComplete` **then**
        **return** `false`
    **for** $r \in P$ **do** ▷ Try every possible pickup
        **if** $r \notin f$ **then**
            $t'' \leftarrow \min\left(t' + t_{r_m,r}, e_r\right)$ ▷ Arrival time at $r$
            **if** $t'' \leq l_r$ and $(q' + q_r) \leq Q$ **then**
                EXTEND($f + \langle r \rangle, t'', q' + q_r$) ▷ Extend with pickup $r$
    **return** `true` ▷ We can complete $f$
$\Phi \leftarrow \emptyset$ ▷ Initialise global set of fragments
**for** $i \in P$ **do**
    EXTEND($\langle i \rangle, e_i, q_i$)

### *4.2. Generating Load Plans*

To determine if there is a feasible loading plan for the PDPTWMS we also use a recursive algorithm to implicitly enumerate all possible loading plans. When tentatively placing a request in a stack during this algorithm we need to consider:

- The stack must have enough spare space.
- A request cannot go in a stack if it contains incompatible requests.
- If we have multiple empty stacks available for a request, we avoid equivalent solutions by considering only the first empty stack.

Implicit enumeration turns out to be very efficient in our computational experiments because the maximum number of requests in each fragment is relatively small.

*CanLoad* is a recursive function which attempts to load a partial fragment $f$ into the stacks which have current contents defined by state $S$. The function repeatedly processes the next request in $f$ by removing it from its stack if it is a delivery or adding it to a stack if it is a pickup. The procedure avoids placing stack incompatible requests in the same stack. If multiple stacks are empty, only the first empty stack encountered is considered.

**Algorithm 2** Generate load plans

```
function CanLoad(f, S)
    if f is empty then
        return true
    r ← first vertex in f
    f′ ← f with r removed
    if r ∈ D then                                  ▷ Remove delivery and check rest of fragment
        s ← stack that contains r − n
        remove r − n from s
        if CanLoad(f′, S) then
            return true
        add r − n to s                             ▷ Restore the state of S
    else
        DoneEmpty ← false
        for s ∈ S do
            if s is empty then                     ▷ Only try the first empty stack
                if not DoneEmpty then
                    add r to s
                    if CanLoad(f′, S) then
                        return true
                    remove r from s                ▷ Restore the state
                    DoneEmpty ← true
            else                                   ▷ Only try this stack if there is space and LIFO applies
                q_s ← current total load of s
                if I_f(r) ∩ s = ∅ and q_s + q_r ≤ Q̂ then
                    add r to s
                    if CanLoad(f′, S) then
                        return true
                    remove r from s                ▷ Restore the state
    return false
```

*4.3. Fragments for the PDPTWRSMS*

For our overall algorithm to work for PDPTWRSMS we only need to replace *CanLoad* with an appropriate stack checking procedure. First we formulate an IP that finds legal loading plans for a fragment, then we discuss some important performance improvements which avoid calling the IP whenever possible.

When considering a candidate loading plan we do not need to check stack capacity is respected every time we load a request (pickup). It is only necessary to check stack capacity at those points in the fragment when a pickup is followed by a delivery. If the stack capacity constraints are respected at these points then they are respected for the whole fragment. Let $\Pi$ denote a set of sets, where each set represents undelivered requests that appear on board every time the vehicle travels from pickups to deliveries. For example, from the fragment $\langle p_1, p_2, d_1, p_3, d_2, d_3\rangle$, $\Pi = \{\{p_1, p_2\}, \{p_2, p_3\}\}$. We make use of $\Pi$ to check that the loads for the sets in $\Pi$ do not exceed the stack capacities.

Now we introduce an IP formulation to check that $f$ can be loaded for the PDPTWRSMS. This IP has a constant objective function of 0, so it simply determines if a feasible solution exists.

Let $\theta$ denote the set of requests (represented by the pickup vertices) that appear in $f$. We denote by $I$ the set of stack incompatible pairs of requests in $\theta$. Let the binary variables $\Upsilon_{ps}$ equal 1 if some part of request $p \in P$ uses stack $s$. Also, let the integer variables $\Lambda_{ps}$ denote the number of units of request $p$ that are placed on stack $s$. The stack assignment problem can be formulated as follows (we omit the constant objective function).

$$\Lambda_{ps} \leq \left(\min\left\{\hat{Q}, q_p\right\}\right)\Upsilon_{ps}, \qquad \forall p \in \theta, \forall s \in S \tag{1}$$

$$\sum_{s \in S} \Lambda_{ps} = q_p, \qquad \forall p \in \theta \tag{2}$$

$$\sum_{p \in \zeta} \Lambda_{ps} \leq \hat{Q}, \qquad \forall s \in S, \forall \zeta \in \Pi \tag{3}$$

$$\Upsilon_{is} + \Upsilon_{js} \leq 1, \qquad \forall (i, j) \in I, \forall s \in S \tag{4}$$

$$\Upsilon_{ps} \in \{0, 1\}, \Lambda_{ps} \in \mathbb{Z}^+ \qquad \forall p \in \theta, \forall s \in S \tag{5}$$

Constraints (1) link $\Lambda_{ps}$ and $\Upsilon_{ps}$. Constraints (2) guarantee that all units of request $p$ will be loaded into stacks. Constraints (3) guarantee that the stack capacity is respected. Constraints (4) ensure that stack incompatible requests will not be loaded on the same stack.

While this IP is relatively efficient for fragments with small numbers of requests and small numbers of stacks, we check for a feasible load plan often during fragment generation. We therefore apply a number of performance improvements, as follows.

If a fragment is PDPTW feasible we can remove any internal sequence of two vertices where the delivery of a request immediately follows the pickup. This is because we know that the maximum vehicle capacity is never exceeded, so there must be enough room across all stacks to load the request. This reduction can be applied repeatedly. For example, consider the following fragment $\langle p_1, p_2, d_1, p_3, p_4, p_5, d_5, d_4, d_2, d_3\rangle$, we can remove the internal sequence $\langle p_4, p_5, d_5, d_4\rangle$. Now we only need to check if the remaining sequence $\langle p_1, p_2, d_1, p_3, d_2, d_3\rangle$ has a feasible loading plan.

We observe that a fragment has a feasible loading plan if no request is ever placed on another. We can check if this is possible by checking if $\sum_{p \in \zeta} \left\lceil q_p/\hat{Q} \right\rceil \leq |S|$ for all $\zeta \in \Pi$.

Clearly, if the maximum request size is no larger than a stack capacity and the fragment is feasible for the PDPTWMS, then it is feasible for the PDPTWRSMS. As the PDPTWMS check is faster, we apply this check first, and only solve the IP (1)–(5) if the fragment is PDPTWMS load infeasible.

In respect of multi-stack loading feasibility, different fragments may share a similar loading configuration if their pickups and deliveries share similar positions and sizes and are sequenced in the same pattern. For example, suppose requests 1, 2 and 3 are the same size as requests 4, 5 and 6 (respectively), $q_1 = q_4$, $q_2 = q_5$ and $q_3 = q_6$, then if $f_1 = \langle p_1, p_2, p_3, d_1, d_2, d_3\rangle$ has a feasible loading plan, so does $f_2 = \langle p_4, p_5, p_6, d_4, d_5, d_6\rangle$. This can be generalised by observing that if $q_1 \geq q_4$, $q_2 \geq q_5$ and $q_3 \geq q_6$ then if $f_1$ has a feasible loading plan then so does $f_2$. Conversely, if $f_2$ does not have a feasible loading plan, then neither does $f_1$.

More formally, we have a fragment $f = \langle r_1, r_2, ..., r_{2m}\rangle$ and it contains pickups $P_f = \langle p_1, p_2, ..., p_m\rangle$ appearing in that order. We define $\gamma_f(r), r \in P_f$ to be the position of $r$ in $P_f$, with $\gamma_f(r) = \gamma_f(r - n), r \in f \cap D$. We define the pickup/delivery sequence $\alpha(f) = \langle \gamma_f(r_1), \gamma_f(r_2), ..., \gamma_f(r_{2m})\rangle$ and we define the load sequence $\beta(f) = \langle q_{p_1}, q_{p_2}, ..., q_{p_m}\rangle$. If we have two fragment $f_1$ and $f_2$ with $\alpha(f_1) = \alpha(f_2)$ and $f_1$ has a feasible load plan and $\beta(f_2) \leq \beta(f_1)$ (where the comparison holds for each term in the sequence) then $f_2$ has a feasible load plan. Conversely if $f_1$ does not have a feasible load plan and $\beta(f_2) \geq \beta(f_1)$ then $f_2$ does not have a feasible load plan. For every pickup/delivery sequence encountered we can cache all load sequences together with the result of checking for a load plan. We can then check this cache first whenever checking a fragment.

The overall algorithm to check the PDPTWRSMS feasibility, for a fragment that is already known to be PDPTW feasible, can be described in the following steps:

1. Remove internal LIFO sequences.
2. Check if the fragment can be loaded with requests never sharing stacks. Return true if this is possible.
3. Check the cache and return cached result if one is found.
4. Check PDPTWMS loading feasibility. If feasible, cache result and return true.
5. Check PDPTWMS loading feasibility with each stack set to have a capacity of $Q$. If infeasible, cache result and return false.
6. Solve stack assignment IP. Cache result and return.

The aim of this algorithm is to minimise the number of stack assignment IPs solved. It is very effective, since 98.8% of the fragments checked do not require an IP to be solved.

### *4.4. Fragment generation speed ups*

A number of other speed ups are also applied. These include:

1. Cache checks apply for the PDPTWMS as well as the PDPTWRSMS.
2. For any pickup/delivery sequence and load sequence, as well as checking the cache directly, we can reverse both sequences and again check the cache.
3. When checking PDPTWMS load plans for a fragment we can remove internal sequences if there is guaranteed to be enough space on one stack to store the request. So if $p_r$ and $d_r$ appear consecutively in a fragment and the load on board immediately before picking up $p_r$ is $q'$, then we can remove $p_r$ and $d_r$ if $q_r + \lfloor q'/k \rfloor \leq \hat{Q}$
4. When checking PDPTWMS load plans for a fragment, if we deliver a request when at most $k$ requests have been picked up, then the request can always be placed on a stack of its own. The pickup and delivery of the request can be removed from the fragment. A similar condition applies for pickups which are followed by at most $k$ deliveries.
5. For the PDPTWMS, when extending a fragment to a new pickup in Algorithm 1 we can first check if the current load on board the vehicle can be packed into the stacks. This packing check makes use of a cache and is very efficient. When the current load cannot be packed, it saves the work of exploring all the possible ways of delivering the current load and determining that none of them have a PDPTWMS load plan.
6. Several heuristic checks can determine that there is no PDPTWRSMS load plan, and so save us having to build and solve the PDPTWRSMS load plan IP. Specifically, for each set of requests that are on board the vehicle we consider the stack incompatibilities and check the following:
    - If all the requests that are incompatible with all other requests have their own stacks, is enough space available for the remaining requests?
    - If a specific request $i$ has its own stack, is there enough space for all requests incompatible with $i$?
    - If two incompatible requests $i$ and $j$ each have their own stack, is there enough space for all the requests incompatible with both $i$ and $j$?

### *4.5. Extended fragments and dominance rules*

We convert the fragments to extended fragments by applying Definition 2.1. Henceforth in this paper we will refer to extended fragments as simply fragments and $\Phi$ will be the set of all extended fragments.

Consider a fragment $f \in \Phi$ that consists of a sequence of pickup and delivery vertices, starting at a pickup vertex and finishing at a pickup vertex or the depot. Let $\sigma_f$ denote the cumulative distances starting from the first vertex and ending at the last vertex of $f$ (the cost of fragment $f$). Let $f^- \in P$ and $f^+ \in P \cup \{0\}$ denote the first vertex and the last vertex of fragment $f$ respectively. Let $\eta_f$ denote the earliest arrival time at $f^+$; $\eta_f$ can be obtained by calculating the arrival time at each vertex working forward from the earliest departure time ($e$) of $f^-$. Also, let $\psi_f$ denote the latest possible departure time from $f^-$; $\psi_f$ can be obtained by calculating the departure time at each vertex working backward from the latest arrival time ($l$) of $f^+$. Finally, if $\theta_{f_1} = \theta_{f_2}$, $f_1^- = f_2^-$, $f_1^+ = f_2^+$, $\sigma_{f_1} \leq \sigma_{f_2}$, $\eta_{f_1} \leq \eta_{f_2}$, and $\psi_{f_1} \geq \psi_{f_2}$, and at least one of the inequalities is strict, then fragment $f_1$ dominates fragment $f_2$. If all the comparisons are equal then $f_1$ dominates $f_2$ if it appears earlier in lexicographical order. All dominated fragments are discarded from $\Phi$.

## 5. Solution approach

Our overall solution approach involves formulating a relaxed model and solving it using the 'Callback' facility of a solver to implement a branch and cut approach. We improve performance by using valid inequalities and incorporating some of the acceleration techniques described in Rist and Forbes (2021) and Rist and Forbes (2022).

### *5.1. Chaining the fragments in a relaxed model*

The models to chain the fragments for PDPTWMS or for PDPTWRSMS are the same. To chain the generated fragments we formulate a network flow model with side constraints. We relax the time windows and the sub-tour (cycle) constraints. We call this model the relaxed PDPTWMS (PDPTWMS-R). Any feasible solution of the PDPTWMS is a feasible solution of the PDPTWMS-R, but some feasible solutions of the PDPTWMS-R may not be feasible solutions of the PDPTWMS as once we chain together fragments it may not be possible to meet all time window constraints or sub-tours may occur. We use lazy constraints to cut off any infeasible solutions generated while solving the resultant integer program.

The main variables of our problem represent fragments and we will also make use of *start* arcs that connect the depot to pickup vertices. Each feasible route for PDPTWMS can be represented as a start arc followed by a series of fragments, the last of which finishes at the depot.

Let $\alpha \in \mathcal{A}$ denote start arcs each of which has a start "tail" vertex (the depot) denoted by $\alpha^- \in V$ and an end "head" vertex denoted by $\alpha^+ \in P$. Likewise a fragment $f \in \Phi$ has a start "tail" vertex $f^- \in P$ and an end "head" vertex $f^+ \in P \cup \{0\}$. Let $c_\alpha$ denote the cost of traversing arc $\alpha \in \mathcal{A}$. Let $\sigma_f$ denotes the cost of fragment $f$. Also, let $\Phi(p)$ denote the set of fragments that contain request $p$ such that $\Phi(p) = \{f \in \Phi, p \in f, p \neq f^+\}, \forall p \in P$. The binary variables $x_f$ are equal to 1 if and only if the fragment $f \in \Phi$ is used in the solution. Finally, the binary variables $y_\alpha$ are equal to 1 if and only if the start arc $\alpha \in \mathcal{A}$ is used in the solution.

The PDPTWMS-R can be formulated as an integer programming problem as follows.

$$\text{minimise} \quad z = \sum_{\alpha \in \mathcal{A}} c_\alpha y_\alpha + \sum_{f \in \Phi} \sigma_f x_f \tag{6}$$

subject to

$$\sum_{f \in \Phi(p)} x_f = 1, \qquad \forall p \in P \tag{7}$$

$$\sum_{\alpha \in \mathcal{A}, \alpha^+ = i} y_\alpha + \sum_{f \in \Phi, f^+ = i} x_f = \sum_{f \in \Phi, f^- = i} x_f, \qquad \forall i \in P \tag{8}$$

$$y_\alpha, x_f \in \{0, 1\}, \qquad \forall \alpha \in \mathcal{A}, \forall f \in \Phi \tag{9}$$

The objective function (6) minimises the total routing cost. Constraints (7) ensure that each request is served exactly once. Finally, constraints (8) ensure that for each pickup vertex in the network the number of start arcs and fragments entering a vertex equals the number of fragments leaving the vertex (flow balance).

We solve this model using a general-purpose integer programming (IP) solver. Using a solver 'Callback' function, we check that every chain of fragments is a valid vehicle route whenever we find a candidate improved integer feasible solution. If all paths in the current IP candidate solution are feasible with respect to time windows we accept the candidate solution. If not, lazy constraints are used to cut off the combination of variables corresponding to the smallest portion of the chain that is infeasible, and to eliminate any cycles that may occur.

We only check routes that consist of two or more fragments, as each route of only one fragment must be feasible for the PDPTWMS. We check each route against the time window constraints and we separately check for cycles. Suppose a route $R$ is infeasible with respect to the time windows condition. Let $f(R)$ represent the set of fragments $f \in \Phi$. We can search the route for the shortest infeasible chain of fragments denoted by $r$. Let $f(r) \subseteq f(R)$ be the set of fragments in the shortest sequence of fragments that violates the time window condition in the route $R$. The lazy constraint to prohibit the infeasible chain is:

$$\sum_{f \in f(r)} x_f \leq |f(r)| - 1 \tag{10}$$

The same lazy constraint can be used to cut off a cycle $r$.

### 5.2. Valid Inequalities

We define a family of Subset-Row (SSR) inequalities, following Rist and Forbes (2021). We define the requests covered by $f$ as $\text{Reqs}(f) = \{p : p \in f \cap P, p \neq f^+\}$. Given three distinct requests $i$, $j$ and $k$, let

$$\mathcal{F}_{ijk} = \{f : f \in \Phi, |\{i, j, k\} \cap \text{Reqs}(f)| \geq 2\}. \tag{11}$$

Since at most one fragment in each set can be used, the following are valid inequalities.

$$\sum_{f \in \mathcal{F}_{ijk}} x_f \leq 1 \qquad \forall\, i, j, k \in P, i, j, k \text{ distinct}$$

Violated inequalities can be easily separated.

### 5.3. Variable fixing

All the computational examples solved in Section 6 have a prioritised objective, first reducing the number of vehicles required and then the total travel cost. We can take advantage of this by solving the linear relaxation of (7) to (9) with the objective value set to minimise the number of vehicles:

$$\text{minimise} \quad z^v = \sum_{\alpha \in \mathcal{A}} c_\alpha y_\alpha. \tag{12}$$

Clearly if $z^v$ is a lower bound on the number of vehicles required, then so too is $\lceil z^v \rceil$. The following proposition holds from Desaulniers et al. (2020).

**Proposition 5.1.** Let $z^*$ be an upper bound on the optimal objective for a minimisation MIP. Suppose the solution to the linear programming relaxation gives a lower bound $z^{\text{LP}}$ and assigns reduced cost $\rho_f > 0$ to binary variable $x_f$. If $z^{\text{LP}} + \rho_f > z^*$ then $x_f = 0$ in every optimal solution to the MIP.

We take advantage of this proposition to tentatively fix to zero variables which will cause the number of vehicles to increase and to tentatively fix to zero those variables which will increase the MIP objective $z^{\text{MIP}}$ such that $z^{\text{MIP}} > (1+\delta)z^{\text{LP}}$ where $\delta$ is a parameter chosen to be 0.025 in our computational testing. This parameter is our "guess" for the size of the IP/LP gap.

In the overall solution approach described below, we first minimise the number of vehicles then minimise the travel cost, with the number of vehicles constrained and assuming that the vehicle cost is 0. When we say "Solve the MIP" we are referring to solving (6)–(9), together with any additional vehicle bounds or variable fixing constraints, using a callback to separate violated constraints (10).

1. Solve the LP relaxation of (7)–(9) with objective (12) to calculate $z^v$.
2. Set $n^v = \lceil z^v \rceil$. $\Phi^v = \{f : f \in \Phi, \rho_f > n^v - z^v\}$. Add the constraint $\sum_{f \in \Phi^v} x_f = 0$.
3. Add a constraint to fix the number of vehicles to $n^v$.
4. Repeatedly solve the LP relaxation of (6)–(9), adding violated SSR inequalities until no more violations are found.
5. If the LP is infeasible then increment $n^v$, remove the constraints added at steps 2 and 3 and return to step 3.
6. Set $\Phi^0 = \{f : f \in \Phi, \rho_f > \delta \cdot z^{\text{LP}}\}$. Add the constraint $\sum_{f \in \Phi^0} x_f = 0$.
7. Solve the MIP.
8. If the MIP is feasible and $z^{\text{MIP}} \leq (1+\delta)z^{\text{LP}}$, stop with the optimal solution.
9. Replace the constraint $\sum_{f \in \Phi^0} x_f = 0$ with $\sum_{f \in \Phi^0} x_f \geq 1$.

10. If the MIP from step 7 was feasible:
    - Set $\Phi^0 = \{f : f \in \Phi, \rho_f > z^{\text{MIP}} - z^{\text{LP}}\}$. Add the constraint $\sum_{f \in \Phi^0} x_f = 0$.
    - Solve the MIP and either improve on the solution from step 7 or show no such improvement is possible. In any event, stop with the optimal solution.
11. Solve the MIP. If it is feasible, stop with the optimal solution. Otherwise, increment $n^v$ remove all extra constraints and return to step 3.

Across all our computational testing, $n^v$ needed to be incremented at most once, and the MIP needed to be solved at most 4 times.

## 6. Computational results

The experiments were run on the Friday cluster at The University of Queensland, on nodes with an Intel Xeon Gold 6240 at 2.60 GHz, supporting AVX-512, running Rocky Linux 9.8. Our methods are implemented in C++ (GCC 12.3.0) and use Gurobi 13.0.0 on a single thread, matching the single-threaded results of Cherkesly and Gschwind (2022). Each instance was given one core, 64GB of memory and a time limit of 3600 seconds, except where noted. All 1,914 runs were executed as a single batch under identical settings.

Of the 1,914 runs, 1,910 were solved to proven optimality and 4 reached the one hour limit; none exhausted the 64GB of memory, and the 4 that exceeded the time limit are right-censored rather than failed, closing when the limit is raised to eight hours as Section 6.1.1 reports. Total solver effort was 9.2 core-hours. Solve times are heavily skewed. The median run takes 0.07 seconds, while only 41 runs exceed one minute and 8 exceed ten minutes. Group means should therefore be read alongside the medians, as a handful of a280 instances carry them.

As is standard for the instances in our tests, the travelling costs $c_{ij}$ are rounded down to the second decimal digit and the travel times $t_{ij}$ are multiplied by 100, rounded up and stored as integers to avoid rounding issues. All time windows are also multiplied by 100 and stored as integers.

The previous state-of-the-art results for solving the PDPTWMS are those reported by Cherkesly and Gschwind (2022). Their results were generated on a computer with an Intel Core i7-6900k processor using CPLEX 12.6.1 to solve the restricted master problem. A time limit of 3600 seconds was imposed for the solution of each instance.

All work to date on the PDPTWMS and its variants uses the instances from Cherkesly et al. (2016). They generated two classes of instances C1 and C2. In the C1 class each request has a size of one unit and the total capacity of each vehicle is six. In the C2 class the size of each request is a random number between three and nine, although the largest request size occurring in the distributed instances is eight. That distinction matters in Section 7.3, because it means every request fits within a single stack in every configuration with up to three stacks. The total capacity of the vehicle is 24 if the number of the vehicle compartments is one or two, and 27 if the vehicle has three compartments. They tested these instances with the number of requests ranging from 25 to 75 in increments of 5. These PDPTWMS instances are classified as follows:

- "a280" instances derived from the a280 instances of the travelling salesman problem library (TSPLIB). There is a total of 198 instances, 99 of each class (C1 and C2). Classified into nine groups, each group consisting of 11 instances depending on the number of requests. For each group the earliest time, at which service at vertex $i \in V \setminus \{0\}$ can begin, is randomly generated in three different settings, that is, 500–1000, 1000–1200, and 1500–2000. For each setting three different time windows lengths were assigned, that is, 15, 30, and 45 as follows.
    1. $e_i \leq 500$, $\forall i \in P$ and $e_i \leq 1000$, $\forall i \in D$, denoted by: w15–500–1000, w30–500–1000 and w45–500–1000.
    2. $e_i \leq 1000$, $\forall i \in P$ and $e_i \leq 1200$, $\forall i \in D$, denoted by: w15–1000–1200, w30–1000–1200 and w45–1000–1200.
    3. $e_i \leq 1500$, $\forall i \in P$ and $e_i \leq 2000$, $\forall i \in D$, denoted by: w15–1500–2000, w30–1500–2000 and w45–1500–2000.

To illustrate, for group w15–500–1000, the length of the time window at each vertex $i \in V \setminus \{0\}$ is 15, and the maximum earliest time, that is randomly generated, is 500 for each pickup vertex and 1000 for each delivery vertex.

- "brd14051", "d18512", "fnl4461", and "nrw1379" instances, 55 instances each, derived from the brd14051, d18512, fnl4461, and nrw1379 instances of the TSPLIB. The total is 440 instances, 220 of each class (C1 and C2). Classified into five groups, each group consisting of 11 instances depending on the number of requests. For each group the earliest time, at which service at vertex $i \in V \setminus \{0\}$ can begin, is randomly generated, and five different time windows lengths were assigned, that is 45, 60, 75, 90, and 120 as follows.
  $e_i \leq 3000,\ \forall i \in P$ and $e_i \leq 4000,\ \forall i \in D$, denoted by: w45–3000–4000, w60–3000–4000, w75–3000–4000, w90–3000–4000, and w120–3000–4000.

Each class was tested with one, two and three stacks. The total number of instances is 1,914 $(3 \times (198 + 440))$.

### *6.1. The PDPTWMS*

In our tests all 1,914 PDPTWMS instances were solved to optimality within the hour except four, all of them a280 instances of class C1 with two or three stacks, three in the group w45–1500–2000 and one in w30–1500–2000. Those four close when the limit is raised to eight hours, as Section 6.1.1 reports, so every instance of the benchmark is solved to optimality by our method. Detailed results for all instances are provided in the supplementary material.

Table 1 and Table 2 report our results for each group of the a280 instances and of the combined brd14051, d18512, fnl4461 and nrw1379 instances, respectively. These two tables report our results alone. A per-group comparison is not possible because Cherkesly and Gschwind (2022) publish instance-level results only for class C2 with two and three stacks, in the supplementary material accompanying that article; class C1 and the one-stack case are reported there only in aggregate, as the six rows reproduced in Table 3. The comparison against the state of the art is therefore made at the aggregate level, in Table 3.

Table 3 summarises the comparison, with the instances grouped by class and number of stacks. In Table 3 both columns are computed on the same basis. The average is over all 319 instances of each group, with an instance not solved to optimality counted at the 3600 second limit, as in their Table 3. In Table 1 and Table 2, which report our results alone, averages are over the instances solved to optimality.

The columns in Table 1–Table 3 are labelled as follows: instances by group (*Group*), the number of instances in each group solved to optimality (*Solved*), and the average CPU time in seconds (*Sec.*). Note that *Sec.* also includes the CPU time used to generate the non-dominated fragments.

Table 1: Results for the a280 groups of C1 and C2 classes

| | Our method | | | | | |
|---|---|---|---|---|---|---|
| | 1 stack | | 2 stacks | | 3 stacks | |
| *Group* | *Solved* | *Sec.* | *Solved* | *Sec.* | *Solved* | *Sec.* |
| C1 class | | | | | | |
| w15-500-1000 | 11 | 0.0 | 11 | 0.3 | 11 | 0.3 |
| w15-1000-1200 | 11 | 0.0 | 11 | 0.7 | 11 | 0.9 |
| w15-1500-2000 | 11 | 0.2 | 11 | 39.5 | 11 | 78.5 |
| w30-500-1000 | 11 | 0.0 | 11 | 0.9 | 11 | 1.7 |
| w30-1000-1200 | 11 | 0.1 | 11 | 21.4 | 11 | 5.9 |
| w30-1500-2000 | 11 | 0.2 | 11 | 49.3 | 10 | 87.3 |
| w45-500-1000 | 11 | 0.1 | 11 | 3.5 | 11 | 10.1 |
| w45-1000-1200 | 11 | 0.1 | 11 | 39.0 | 11 | 41.8 |
| w45-1500-2000 | 11 | 0.4 | 10 | 568.5 | 9 | 574.3 |
| C2 class | | | | | | |
| w15-500-1000 | 11 | 0.0 | 11 | 0.1 | 11 | 0.2 |
| w15-1000-1200 | 11 | 0.0 | 11 | 1.4 | 11 | 1.0 |
| w15-1500-2000 | 11 | 0.1 | 11 | 7.2 | 11 | 26.4 |
| w30-500-1000 | 11 | 0.0 | 11 | 0.9 | 11 | 0.6 |
| w30-1000-1200 | 11 | 0.1 | 11 | 1.0 | 11 | 3.2 |
| w30-1500-2000 | 11 | 0.2 | 11 | 8.7 | 11 | 58.2 |
| w45-500-1000 | 11 | 0.1 | 11 | 2.0 | 11 | 2.9 |
| w45-1000-1200 | 11 | 0.1 | 11 | 9.2 | 11 | 20.2 |
| w45-1500-2000 | 11 | 0.2 | 11 | 52.8 | 11 | 103.9 |

Table 2: Results for the brd14051, d18512, fnl4461, and nrw1379 groups of C1 and C2 classes

| | Our method | | | | | |
|---|---|---|---|---|---|---|
| | 1 stack | | 2 stacks | | 3 stacks | |
| *Group* | *Solved* | *Sec.* | *Solved* | *Sec.* | *Solved* | *Sec.* |
| C1 class | | | | | | |
| w45-3000-4000 | 44 | 0.0 | 44 | 0.2 | 44 | 0.1 |
| w60-3000-4000 | 44 | 0.0 | 44 | 0.2 | 44 | 0.2 |
| w75-3000-4000 | 44 | 0.0 | 44 | 0.5 | 44 | 0.5 |
| w90-3000-4000 | 44 | 0.0 | 44 | 0.6 | 44 | 0.7 |
| w120-3000-4000 | 44 | 0.0 | 44 | 0.8 | 44 | 1.8 |
| C2 class | | | | | | |
| w45-3000-4000 | 44 | 0.0 | 44 | 0.1 | 44 | 0.2 |
| w60-3000-4000 | 44 | 0.0 | 44 | 0.3 | 44 | 0.3 |
| w75-3000-4000 | 44 | 0.0 | 44 | 0.3 | 44 | 0.4 |
| w90-3000-4000 | 44 | 0.0 | 44 | 0.4 | 44 | 0.8 |
| w120-3000-4000 | 44 | 0.0 | 44 | 0.9 | 44 | 1.5 |

Table 3: Comparison summary results

| | Cherkesly and Gschwind (2022) | | Our method | |
|---|---|---|---|---|
| *Group* | *Solved* | *Sec.* | *Solved* | *Sec.* |
| C1 class | | | | |
| 1 stack | 317 | 56.6 | 319 | 0.1 |
| 2 stacks | 263 | 809.1 | 318 | 34.8 |
| 3 stacks | 258 | 864.2 | 316 | 58.1 |
| C2 class | | | | |
| 1 stack | 319 | 12.7 | 319 | 0.1 |
| 2 stacks | 288 | 488.9 | 319 | 3.2 |
| 3 stacks | 286 | 503.0 | 319 | 7.9 |
| Total | 1,731 | | 1,910 | |
| Newly solved | | | | 179 |

The results for the PDPTWMS instances indicate that our method outperforms the branch-price-and-cut algorithm of Cherkesly and Gschwind (2022). It solves 1,910 of the 1,914 instances to optimality within the hour against 1,731, closing at least 179 that were previously unsolved, and the remaining four with a limit of eight hours. It is also faster in every group of Table 3, though not uniformly, by about 900 times on the single-stack instances of class C1 and 250 times on those of class C2, and by 15 to 155 times with two or three stacks, which is where all four of the instances unsolved within the hour lie. Over all 1,914 instances the mean solution time falls by a factor of about 26.

*6.1.1. The instances not solved within the hour*

Four of the 1,914 runs reached the one hour limit, and all four are class C1 a280 instances of the widest time windows and the latest release times, namely 131–w45–1500–2000 with two stacks and 141–w30–1500–2000, 141–w45–1500–2000 and 151–w45–1500–2000 with three. They are the largest instances of the family at their stack count, and every one of them ran out of time inside the integer program rather than during fragment generation, with a feasible solution already in hand. Only one of the four, 141–w30–1500–2000, had reached the last of its integer programs when the hour expired, and was then within 0.04% of proving optimality. The other three were still solving the first, so how close they stood at that moment says little, since the program that eventually closed them carries the violated inequalities of Section 5.1 and a column set reduced by fixing. On all four a feasible solution with minimum fleet size had been found. What the hour left open on these instances is the routing distance alone. To establish whether the limit rather than the method is what leaves that open, we ran these four again with a limit of eight hours, under the settings of Section 6 in every other respect.

All four then solved to proven optimality, and Table 4 reports them. The times are not extreme. The longest of the four took four hours and thirty-nine minutes and the shortest one hour and twenty minutes, so these instances sit just beyond the one-hour limit rather than far outside the reach of the method. On 141–w30–1500–2000 the solution found within the hour was already optimal and only the proof was missing, which is what a gap of 0.04% suggests; on the other three the longer run improved the incumbent by 0.22%, 3.01% and 3.31%.

Table 4: The four instances not closed within the hour, run again with a limit of eight hours

| *Instance* (a280-C1) | *Stacks* | *Sec.* | *Vehicles* | *Distance* | $F_N$ |
|---|---|---|---|---|---|
| 131–w45–1500–2000 | 2 | 7,663.6 | 9 | 9,524.45 | 3,338,597 |
| 141–w30–1500–2000 | 3 | 4,803.8 | 9 | 11,375.44 | 1,765,531 |
| 141–w45–1500–2000 | 3 | 15,732.6 | 8 | 10,639.63 | 8,853,465 |
| 151–w45–1500–2000 | 3 | 16,738.8 | 9 | 10,826.18 | 14,982,594 |

The fragment counts explain the cost. These four generate between 1.8 and 15.0 million non-dominated fragments, against a median of 11,337 over the a280 instances at two stacks, and both the enumeration and the integer program scale with that number. Every instance of the benchmark is therefore solved to optimality by our method, four of them needing more than the hour that the comparison in Table 3 allows. The counts reported elsewhere in this section keep to the one-hour limit, so that they remain comparable with Cherkesly and Gschwind (2022).

### *6.2. The PDPTWRSMS*

The goal of our testing in this section is to compare the PDPTWRSMS to the PDPTWMS. Also, we show how the techniques described in Section 4.3 help reduce the fragment generation time and the number of stack assignment IPs solved (model (1)–(5)).

We have tested the algorithm of the PDPTWRSMS using the a280 instances. The a280 group is chosen deliberately, because it is much the hardest of the five families and so the most informative test of the techniques of Section 4.3. Its coordinates come from a printed-circuit-board drilling pattern rather than from geography, and the resulting regular layout produces many tied distances, which in turn defeats the dominance rules and inflates the fragment set. At two stacks the median a280 instance yields 11,337 non-dominated fragments against 1,818 for the other four families, a sixfold difference, and the median solution time is 0.92 s against 0.09 s. Every one of the instances not closed within the hour in Section 6.1.1 is an a280 instance. The other four families are drawn from towns in Germany and are geometrically irregular; they are solved so easily that they would not exercise the fragment-generation speed-ups at all. Because the number of instances is large, we limit our test to instances of 70 and 75 requests (141 and 151 vertices). We only considered instances of C2 class because in the instances of C1 class each request has a size of one unit.

The columns of the tables in this section represent the following: *Instance*: the instance name; $F_N$: the total number of non-dominated fragments; *IPs*: the number of stack assignment IPs solved; $\text{Sec}_{\text{IP}}$: the CPU time in seconds used to solve the stack assignment IPs, rounded to the nearest tenth; $\text{Sec}_f$: the CPU time in seconds used to generate non-dominated fragments, including $\text{Sec}_{\text{IP}}$ time, rounded to the nearest tenth; *Sec.*: the total CPU time in seconds used to solve the instance, including $\text{Sec}_f$ time, rounded to the nearest tenth; $Z^*$: the optimal solution obtained within the time limit or left blank if no feasible solution was found. If only a feasible solution was obtained within the time limit, then it is reported in bold in column $Z^*$ as a best upper bound found; *Gap*: the percentage gap computed as $100(Z^* - Z_R)/Z^*$, where $Z_R$ is the lower bound at the root node.

#### *6.2.1. The added value of stack splitting*

Table 5 and Table A.1 provide computational results for the instances of two and three stacks respectively, both at a vehicle capacity of 24. The benchmark gives a three-stack vehicle a capacity of 27, which we keep for the PDPTWMS results of Section 6 so that the comparison with Cherkesly and Gschwind (2022) is on the published instances; here we hold the capacity at 24 instead, so that the two configurations differ in the number of stacks alone. The columns $S_{\text{Veh}}$ and $S_{\text{Dist}}$ show the saving against the PDPTWMS in the number of vehicles and in the total travelled distance respectively. If the number of vehicles used for the PDPTWMS and the PDPTWRSMS are $\text{Veh}_{\text{MS}}$ and $\text{Veh}_{\text{MS-S}}$, then $S_{\text{Veh}} = \text{Veh}_{\text{MS-S}} - \text{Veh}_{\text{MS}}$, a count rather than a percentage, since splitting never saves more than one vehicle on any instance here and expressing that as a percentage of a fleet whose size varies by instance obscures it. The distance saving is a percentage, $S_{\text{Dist}} = 100(\text{Dist}_{\text{MS-S}} - \text{Dist}_{\text{MS}})/\text{Dist}_{\text{MS}}$.

Analysing the results in Table 5 and Table A.1, it can be seen that a reduction in cost was obtained for all instances. The distance increases for some instances when the number of vehicles is reduced. This is because the vehicle cost is relatively high, and the priority is to minimise the number of vehicles. Moreover, we see that a reduction of the number of vehicles was more common in the three stack case (13 out of 18 instances) than in the two stack case (3 out of 18 instances). This suggests that splitting customer requests is more effective as the number of stacks increases from 2 to 3, which may be because the stack capacity has reduced – in this case from 12 to 8.

The stack splitting extension usually leads to increased solution time. This is due to the flexibility in the loading configuration, thus more fragments are generated. For example, in the most difficult instance – 151-w45-1500-2000 with 3 stacks – the PDPTWMS had 2,328,776 non-dominated fragments whereas the PDPTWRSMS had 3,617,549.

Table 5: PDPTWRSMS for 2 stacks vehicle with capacity 24

| *Instance* (a280-C2) | $F_N$ | *IPs* | *Gap* | $\text{Sec}_{\text{IP}}$ | $\text{Sec}_f$ | *Sec.* | $Z^*$ | $S_{\text{Veh}}$ | $S_{\text{Dist}}$ |
|---|---|---|---|---|---|---|---|---|---|
| 141-w15-500-1000 | 4,265 | 78 | 5.285 | 0.0 | 0.1 | 0.9 | 2,013,835.44 | 0 | −2.01 |
| 151-w15-500-1000 | 3,197 | 35 | 0.001 | 0.0 | 0.0 | 0.2 | 2,116,833.61 | −1 | 7.89 |
| 141-w15-1000-1200 | 8,809 | 15 | 3.176 | 0.0 | 0.1 | 1.0 | 1,813,158.94 | −1 | −3.96 |
| 151-w15-1000-1200 | 12,387 | 119 | 3.085 | 0.0 | 0.1 | 3.6 | 1,714,685.07 | 0 | −6.13 |
| 141-w15-1500-2000 | 64,511 | 294 | 3.836 | 0.1 | 0.7 | 23.7 | 1,312,428.75 | 0 | −2.07 |
| 151-w15-1500-2000 | 164,856 | 1,427 | 5.923 | 0.3 | 2.1 | 22.6 | 1,413,100.98 | 0 | −4.55 |
| 141-w30-500-1000 | 11,277 | 141 | 3.300 | 0.0 | 0.1 | 1.9 | 1,713,208.73 | 0 | −0.64 |
| 151-w30-500-1000 | 13,731 | 141 | 4.119 | 0.0 | 0.1 | 3.8 | 1,713,337.32 | 0 | −0.31 |
| 141-w30-1000-1200 | 30,372 | 155 | 4.705 | 0.0 | 0.3 | 3.5 | 1,512,482.58 | 0 | −3.36 |
| 151-w30-1000-1200 | 28,306 | 229 | 4.709 | 0.0 | 0.3 | 3.2 | 1,612,795.20 | 0 | −3.55 |
| 141-w30-1500-2000 | 348,019 | 2,611 | 3.006 | 0.5 | 5.2 | 78.4 | 1,212,488.76 | 0 | −9.59 |
| 151-w30-1500-2000 | 362,711 | 1,553 | 5.660 | 0.3 | 5.1 | 33.1 | 1,312,367.22 | 0 | −5.34 |
| 141-w45-500-1000 | 15,014 | 205 | 1.867 | 0.0 | 0.2 | 2.8 | 1,713,586.67 | 0 | −1.56 |
| 151-w45-500-1000 | 43,435 | 567 | 3.187 | 0.1 | 0.6 | 15.9 | 1,713,457.57 | 0 | −0.59 |
| 141-w45-1000-1200 | 54,305 | 659 | 2.420 | 0.1 | 0.8 | 50.0 | 1,312,423.73 | −1 | 6.06 |
| 151-w45-1000-1200 | 120,321 | 626 | 6.688 | 0.1 | 1.6 | 16.1 | 1,412,202.80 | 0 | −0.44 |
| 141-w45-1500-2000 | 408,835 | 4,566 | 7.341 | 0.8 | 8.3 | 250.0 | 1,311,960.88 | 0 | −0.76 |
| 151-w45-1500-2000 | 1,529,017 | 12,309 | 3.814 | 2.9 | 25.1 | 378.0 | 1,212,236.30 | 0 | −9.02 |

#### *6.2.2. The impact of increasing the number of stacks with stack splitting*

Table 6 presents the impact that increasing the number of stacks from two to three has on the number of vehicles, the total travelled distance, the number of generated fragments, and the total computational time.

The columns in Table 6 represent the following: $\Delta_{F_N}$: the percentage difference in the number of non-dominated fragments, computed as $100\left(F_{N_3} - F_{N_2}\right)/F_{N_2}$, where $F_{N_2}$ and $F_{N_3}$ are the number of non-dominated fragments with two and three stacks respectively; $\Delta_{\text{Sec}}$: the percentage difference in the total CPU time in seconds used to solve the instance, computed as $100(\text{Sec}_3 - \text{Sec}_2)/\text{Sec}_2$, where $\text{Sec}_2$ and $\text{Sec}_3$ are the total CPU time in seconds used to solve the instance with two and three stacks respectively; $\Delta_{\text{Dist}}$: the percentage difference in the total travelled distance, computed as $100(\text{Dist}_3 - \text{Dist}_2)/\text{Dist}_2$, where $\text{Dist}_2$ and $\text{Dist}_3$ are the total travelled distance with two and three stacks respectively; *Vehicles*: the difference in the number of vehicles as a count, $\text{Veh}_3 - \text{Veh}_2$, which $\Delta_{\text{Veh}}$ alone cannot be read back to because the fleet size differs by instance; $\Delta_{\text{Veh}}$: the percentage difference in the number of vehicles, computed as $100(\text{Veh}_3 - \text{Veh}_2)/\text{Veh}_2$, where $\text{Veh}_2$ and $\text{Veh}_3$ are the number of vehicles with two and three stacks respectively.

Table 6 shows the positive impact on cost of increasing the number of stacks. In particular, the number of vehicles has decreased in nearly every instance. However, due to the loading flexibility we observe a sharp increase in the number of fragments. This makes the problem harder to solve and usually increases the solution time. The time to generate fragments $\text{Sec}_f$ is higher for all instances.

One of the PDPTWMS assumptions is that a request's size may not exceed the stack capacity. Dividing the vehicle capacity of 24 into four stacks breaks that assumption, and without splitting no instance then has a solution. The stack depth falls to six, so the requests of size seven and eight, which every class C2 instance contains, fit in no stack. Section 7.3 reports the PDPTWRSMS results for this configuration over all 319 class C2 instances, and Table 10 compares them with two and three stacks at the same capacity.

Table 6: Two stacks against three stacks, both at capacity 24

| *Instance* (a280-C2) | $\Delta_{F_N}$ | $\Delta_{\text{Sec}}$ | $\Delta_{\text{Dist}}$ | $\Delta_{\text{Veh}}$ | *Vehicles* |
|---|---|---|---|---|---|
| 141-w15-500-1000 | 28.8 | −65.9 | −4.8 | −5.0 | −1 |
| 151-w15-500-1000 | 40.4 | 595.8 | −11.7 | −4.8 | −1 |
| 141-w15-1000-1200 | 53.1 | −8.7 | −3.6 | −5.6 | −1 |
| 151-w15-1000-1200 | 36.0 | −45.2 | −3.8 | −5.9 | −1 |
| 141-w15-1500-2000 | 130.7 | 157.2 | −7.3 | −7.7 | −1 |
| 151-w15-1500-2000 | 116.6 | 524.8 | −1.1 | −14.3 | −2 |
| 141-w30-500-1000 | 35.1 | 137.1 | −4.4 | −5.9 | −1 |
| 151-w30-500-1000 | 42.2 | −8.6 | −3.5 | −5.9 | −1 |
| 141-w30-1000-1200 | 75.8 | 44.4 | −8.0 | −6.7 | −1 |
| 151-w30-1000-1200 | 79.2 | 685.4 | −4.9 | −6.2 | −1 |
| 141-w30-1500-2000 | 136.1 | 623.3 | −12.3 | 0.0 | 0 |
| 151-w30-1500-2000 | 118.0 | 291.3 | −2.0 | −7.7 | −1 |
| 141-w45-500-1000 | 75.0 | 143.5 | −4.1 | −5.9 | −1 |
| 151-w45-500-1000 | 52.7 | 90.1 | −7.2 | −5.9 | −1 |
| 141-w45-1000-1200 | 61.4 | −30.7 | −10.5 | 0.0 | 0 |
| 151-w45-1000-1200 | 65.7 | 244.4 | 1.3 | −7.1 | −1 |
| 141-w45-1500-2000 | 182.1 | 65.4 | −4.9 | −7.7 | −1 |
| 151-w45-1500-2000 | 136.6 | 716.5 | −8.7 | 0.0 | 0 |

#### *6.2.3. The effectiveness of the techniques described in Section 4.3*

Recall that the overall algorithm to find a PDPTWRSMS loading plan uses the following steps:

1. Remove internal LIFO sequences.
2. Check if the fragment can be loaded with requests never sharing stacks. Return true if this is possible.
3. Check the cache and return cached result if one is found.
4. Check PDPTWMS loading feasibility. If feasible, cache result and return true.
5. Check PDPTWMS loading feasibility with each stack set to have a capacity of $Q$. If infeasible, cache result and return false.
6. Solve stack assignment IP. Cache result and return.

Each step is a sufficient condition for a loading plan to exist, or a necessary one for it not to, and the last is the complete test the others stand in front of. A call is resolved by the first test that applies, so each call is counted once in Table 7, against whichever test returned the answer. A speed up of Section 4.4 that sharpens a step is counted under that step. Step 5 and the stack assignment IP are reached only when splitting is allowed; without it, step 4 is the complete test and every call ends at or before it.

Table 7: Share of the calls each test is the first to resolve, over the two-stack and three-stack configurations at capacity 24, each solved with and without splitting. Steps are those of the list above, each with the speed ups of Section 4.4 that sharpen it

| *Test* | *PDPTWMS* | *PDPTWRSMS* |
|---|---|---|
| 2 Never share a stack | 31.7% | 29.2% |
| 3 Cache | 63.7% | 67.3% |
| 4 PDPTWMS feasibility | 4.6% | 1.0% |
| – Incompatibility heuristics | – | 1.3% |
| 5 Capacity $Q$ relaxation | – | $< 0.1\%$ |
| 6 Stack assignment IP | – | 1.2% |
| Calls | 27,191,477 | 35,850,642 |

In order to see the effectiveness of these techniques we use the data of Table 5 and Table A.1 for the two-stack and three-stack cases, both at capacity 24. In these tables we see that the number of stack assignment IPs solved is much smaller than the number of non-dominated fragments generated. Across all the two-stack instances the number of IPs solved is 0.80% of the number of non-dominated fragments, and for three stacks it is 2.18%. While it is not shown in the tables, our tests also show that the stack assignment IP is infeasible (i.e. no loading plan exists) 68% of the time it is used in two-stack instances and 28% of the time it is used in three-stack instances. The third stack thus reaches the IP far more often, which is how increasing the number of stacks increases the prevalence of splitting; that a smaller share of those calls is infeasible says the same thing, since with shallower stacks more of the fragments that need splitting can actually be split.

We also ran an experiment where we disabled the check for Load Equivalent fragments. These results are shown in Table 8. In this table, column $\Delta_{\text{IP}}$ shows the percentage reduction in the number of stack assignment IPs solved when using the check for Load Equivalent fragments. It is computed as $100(\text{IP}_W - \text{IP}_O)/\text{IP}_O$, where $\text{IP}_O$ is the number of IPs without the check, and $\text{IP}_W$ is the number of IPs with the check. Likewise, column $\Delta_{\text{Sec}_f}$ shows the percentage reduction in the fragment generation time. It is computed as $100\Big(\text{Sec}_{f_W} - \text{Sec}_{f_O}\Big)/\text{Sec}_{f_O}$, where $\text{Sec}_{f_O}$ is the CPU time in seconds used to generate non-dominated fragments without the check, and $\text{Sec}_{f_W}$ is the CPU time in seconds used to generate non-dominated fragments with the check. Checking for Load Equivalent fragments provides a significant speedup, both in reducing the number of stack assignment IPs solved and in reducing fragment generation time. Averaged over the 18 instances the check removes 69.8% of the IPs with two stacks and 78.9% with three, and 38.9% and 57.7% of the fragment generation time. The effect grows with the difficulty of the instance, and is particularly significant for the more difficult problems towards the bottom of the table. Over the four rows at the bottom the check removes 82.4% and 89.9% of the IPs, against 47.6% and 62.1% over the four at the top. It is also worth more with three stacks than with two. The saving in fragment generation time is larger with three on every one of the 18 instances, and the reduction in the number of IPs is larger on 17 of them. That is what the configuration would suggest. The third stack reaches the stack assignment IP far more often, and its IPs are harder, since a request may be divided between more compartments. With the check in place a stack assignment IP takes 0.43 ms on average with three stacks against 0.21 ms with two. Each call the check removes is therefore worth more.

Table 8: The reduction in the IPs and the fragment generation time

| | 2 stacks (capacity 24) | | 3 stacks (capacity 24) | |
|---|---|---|---|---|
| *Instance* (a280-C2) | $\Delta_{\text{IP}}$ | $\Delta_{\text{Sec}_f}$ | $\Delta_{\text{IP}}$ | $\Delta_{\text{Sec}_f}$ |
| 141-w15-500-1000 | −44.7 | −20.9 | −51.4 | −36.8 |
| 151-w15-500-1000 | −53.9 | −25.0 | −50.0 | −40.0 |
| 141-w15-1000-1200 | −25.0 | −8.5 | −75.5 | −30.1 |
| 151-w15-1000-1200 | −66.6 | −31.1 | −71.5 | −39.8 |
| 141-w15-1500-2000 | −84.5 | −41.9 | −88.0 | −58.9 |
| 151-w15-1500-2000 | −83.5 | −53.8 | −89.3 | −71.2 |
| 141-w30-500-1000 | −65.1 | −32.6 | −70.7 | −54.6 |
| 151-w30-500-1000 | −60.2 | −30.6 | −65.9 | −45.9 |
| 141-w30-1000-1200 | −69.2 | −34.3 | −88.5 | −58.5 |
| 151-w30-1000-1200 | −65.8 | −28.5 | −78.2 | −48.9 |
| 141-w30-1500-2000 | −87.1 | −59.7 | −89.7 | −76.5 |
| 151-w30-1500-2000 | −91.3 | −51.8 | −94.1 | −72.4 |
| 141-w45-500-1000 | −55.3 | −32.4 | −69.1 | −47.9 |
| 151-w45-500-1000 | −74.0 | −37.7 | −79.2 | −61.4 |
| 141-w45-1000-1200 | −83.6 | −57.3 | −88.1 | −72.1 |
| 151-w45-1000-1200 | −76.8 | −34.7 | −91.5 | −67.5 |
| 141-w45-1500-2000 | −79.9 | −56.1 | −89.6 | −75.3 |
| 151-w45-1500-2000 | −89.4 | −63.2 | −90.3 | −81.2 |

## 7. Analysis of the optimal solutions

Because our method reports the selected fragments for every instance it closes, the optimal solutions themselves can be examined rather than only their cost. This section does so for class C2 at a total vehicle capacity of 24, and Section 7.5 alone returns to the whole benchmark. Two things narrow it. Splitting a request is the freedom the PDPTWRSMS adds, and it cannot apply to class C1, whose requests are of unit size and cannot be divided, so every comparison that involves splitting is a class C2 comparison. The benchmark also raises the total capacity from 24 to 27 when the third stack is added, so holding the capacity at 24 throughout is what separates the effect of the stack from the effect of the capacity granted alongside it, and the three-stack results below are therefore the capacity-24 runs of Section 6 rather than the benchmark's own. All 319 class C2 instances are solved to optimality at all three stack counts, and their solutions comprise 15,762 vehicle routes. This section asks what those solutions look like — how large a fleet is needed, how much each vehicle carries, what the freedom to split a request across stacks is worth, how that compares with loosening the loading policy instead, and how the routes are assembled from fragments. The last of these bears directly on the separation used in Section 5.1. The questions about splitting and about the loading policy each require solving the instances under configurations other than those of Section 6, and the additional runs are described where they are used.

### *7.1. Fleet size*

Going from one stack to two removes vehicles in every one of the 319 instances, reducing the mean fleet from 21.5 to 14.3, a reduction of 7.2 vehicles on average. The third stack is worth far less and, at a capacity held fixed, it is not always a gain. It reduces the fleet on 56.4% of instances, leaves it unchanged on 39.5% and increases it on 4.1%, for a mean change of −0.639 vehicles. Three stacks of eight are shallower than two of twelve, and the example in Section 3 shows how a loading that two stacks allow can be lost when the same capacity is divided into three.

A single instance shows the effect. In a280–51–c12–w45–1500–2000 the same 25 requests and the same geometry are served by seven vehicles when the vehicle has one stack and by five when it has two. The third stack removes no further vehicle here, but it does shorten the routes, from 4,155.12 to 3,560.11, a saving of 14.3%.

The reason a single stack is so restrictive is a statement about order rather than about capacity, and Figure 2 shows it directly. Each route is drawn as the path traced by the number of requests on board, rising by one at every pickup and falling by one at every delivery, returning to zero when the vehicle empties. The path alone records only how many requests are held. What decides single-stack feasibility is which request each descent removes. Under LIFO loading a vehicle can unload only the item currently on top, so every delivery must close the most recently opened pickup still on board. A route serving $p_a, p_b, d_b, d_a$ satisfies this; $p_a, p_b, d_a, d_b$ does not, because $a$ lies buried beneath $b$, and no reordering of that route repairs it.

The two columns differ qualitatively, and in the same way on both instances. In the one-stack solutions every delivery closes the item on top, because no other is available, and the paths stay correspondingly low. Every route of the two-stack solutions buries a request at some point, and each is therefore infeasible for a single-stack vehicle however its stops are reordered.

The figure also makes clear that this is a restriction of order and not of size. A class C2 vehicle holds 24 units against demands of three to eight, so it has room for between three and eight requests at once. Yet across the one-stack solutions a vehicle carries a mean peak of only 1.93 requests, well short of what it could physically hold. The binding constraint at one stack is LIFO, not capacity. That is what forces the one-stack fleet on that instance to seven vehicles. The vehicle is not too small; the requests worth combining cannot be sequenced to satisfy LIFO on a single stack.

### *7.2. What a vehicle carries*

The mean number of requests carried by a vehicle rises from 2.33 with one stack to 3.50 with two and 3.67 with three, and the largest route grows from seven requests to ten. Those figures describe the whole route; what changes underneath them is how much the vehicle holds at any one moment, and that is where the effect of the stacks is visible.

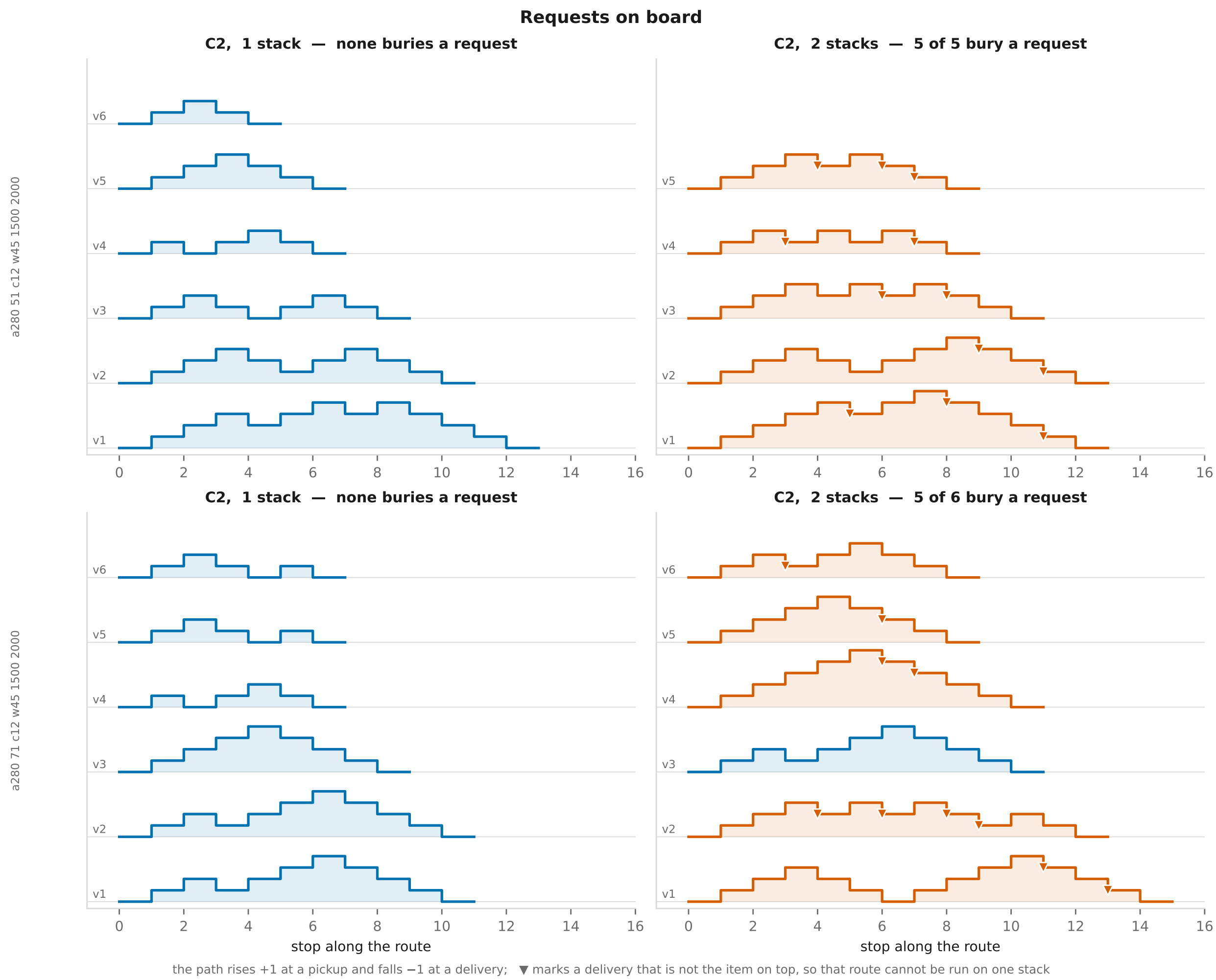

Figure 2: Requests on board along each route, for a280–51–c12–w45–1500–2000 (top) and a280–71–c12–w45–1500–2000 (bottom). The path rises +1 at a pickup and falls −1 at a delivery; ▼ marks a delivery of an item that is not on top, which one stack cannot perform.

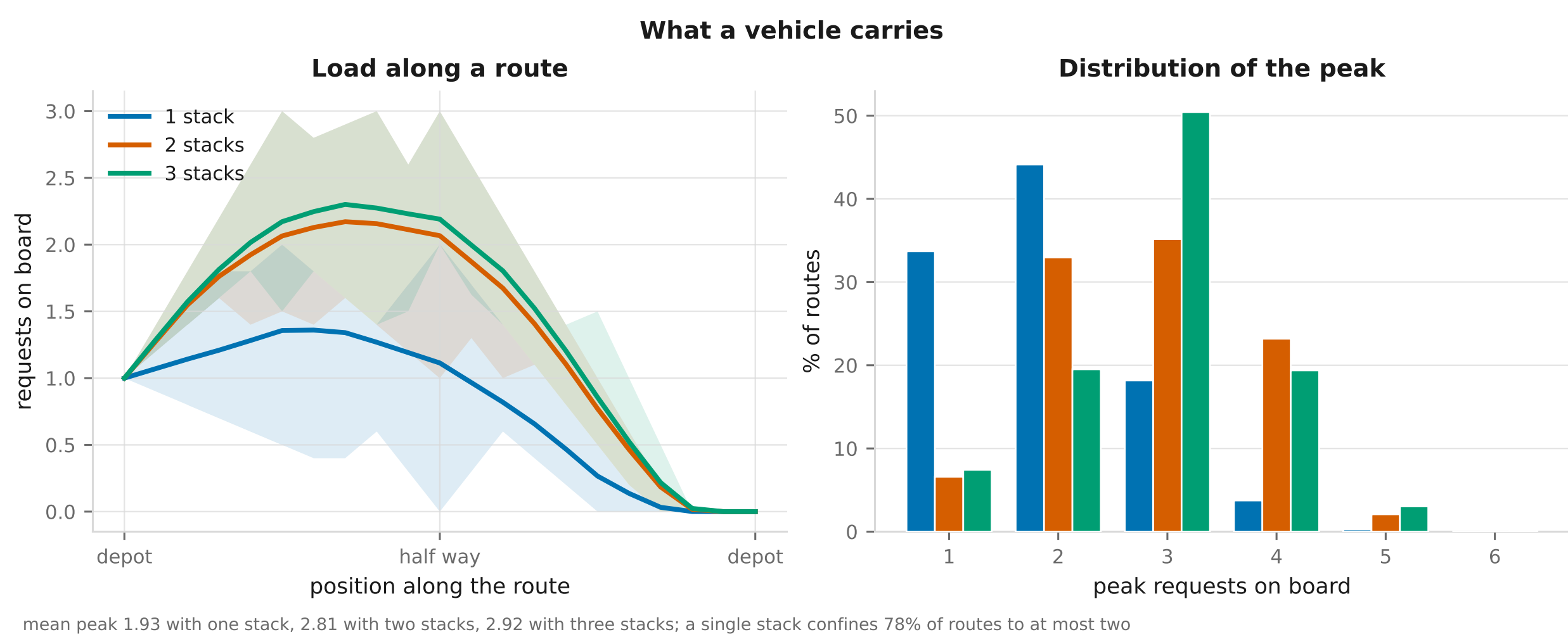

Figure 3: Requests on board along a route, averaged over all routes with the interquartile band (left), and the distribution of each route's peak (right).

The left panel of Figure 3 traces the load along a route. With a single stack it barely rises. The mean peaks at 1.36 requests about a third of the way along the route and the curve is almost flat, because LIFO loading forces the vehicle to deliver an item before it can collect anything that would sit on top. The vehicle spends most of its route carrying one request, so a route is close to a sequence of individual pickup and delivery pairs. With two stacks the same curve peaks at 2.17, about 60% higher, and has a pronounced arc, with the vehicle accumulating requests over the first half of its route and unwinding them over the second. That arc is the consolidation which removes vehicles from the fleet.

The right panel shows the peak each route reaches. A single stack confines 77.8% of routes to at most two open requests, and only 4.0% ever hold four or more. With two stacks the mode moves from two to three, 25.3% of routes reach four or more, and 2.1% reach five or six. The mean peak rises from 1.93 to 2.81.

The three-stack distribution sits a little above the two-stack one, at a mean peak of 2.92 against 2.81, so the third stack does let a vehicle hold more at once, but the step it makes is a fraction of the one the second stack makes. Fewer routes reach four open requests or more, 22.6% against 25.3%, because three stacks of eight bury less than two of twelve but hold less in each as well.

### *7.3. The value of request splitting*

Splitting a request across stacks is the freedom the PDPTWRSMS adds, so it is worth asking what it buys. We solved the 319 instances with splitting enabled in each of the three multi-stack configurations at capacity 24, and compare against the same instances without it wherever a solution without it exists. A single-stack vehicle has nothing to divide a request across, so the comparison begins at two stacks.

It is worth being clear about what splitting is not doing here. The largest request in the benchmark is of size eight, and the per-stack capacity is twelve with two stacks and eight with three at capacity 24. Every request therefore fits within a single stack in each of these configurations, and splitting is never required for feasibility. What it provides is freedom in packing, since a request may be divided to make room rather than because it would not otherwise fit.

The consequence is that the gains reported in Table 9 are a lower bound on what splitting is worth. In these configurations it can only improve a packing that was already feasible; it is never the difference between a solution and none. Requests larger than a stack — the case motivating the variant — arise from the same instances once the vehicle is divided into four stacks, and the end of this section takes that case up.

Because splitting only enlarges the set of feasible loadings, it cannot make a solution worse, and the computations bear this out. Over the 319 instances solved in all four runs, two stacks and three with and without splitting, splitting increases the fleet on none of them. Table 9 reports what it gains.

Table 9: Effect of request splitting, over the 319 class C2 instances solved at two stacks and at three, each with and without splitting

| *Configuration* | *Stack depth* | *Mean fleet change* | *Saves a vehicle* | *Cost at equal fleet* |
|---|---|---|---|---|
| 2 stacks, capacity 24 | 12 | −0.154 | 15.0% | −1.56% |
| 3 stacks, capacity 24 | 8 | −0.323 | 30.1% | −2.27% |

The gain grows as the stacks get shallower, and this is the natural reading, since a shallow stack is what makes a request awkward to place, so the freedom to divide one is worth most exactly where the stacks are shallowest. With two stacks, each half the vehicle, splitting saves a vehicle on 15.0% of instances; with three stacks of a third of the vehicle each, on 30.1%. Where it removes a vehicle it almost always removes exactly one, two going on a single instance with two stacks and on seven with three.

Splitting and additional stacks are alternative responses to the same difficulty, and they can be compared from a common baseline of two stacks at capacity 24. Figure 4 holds the capacity at 24 throughout, and its top right panel makes that comparison.

At constant capacity both changes help, and they compose. Allowing requests to be split reduces the mean fleet by 0.154 vehicles and adding a third stack by 0.639; together they reduce it by 0.962, more than either alone and more than their sum. A fourth stack cannot be added without splitting at all, and with it the fleet falls by 0.878 from the same baseline — a little less than the third stack with splitting achieves. The end of this section examines why the fourth stack gives some of the gain back.

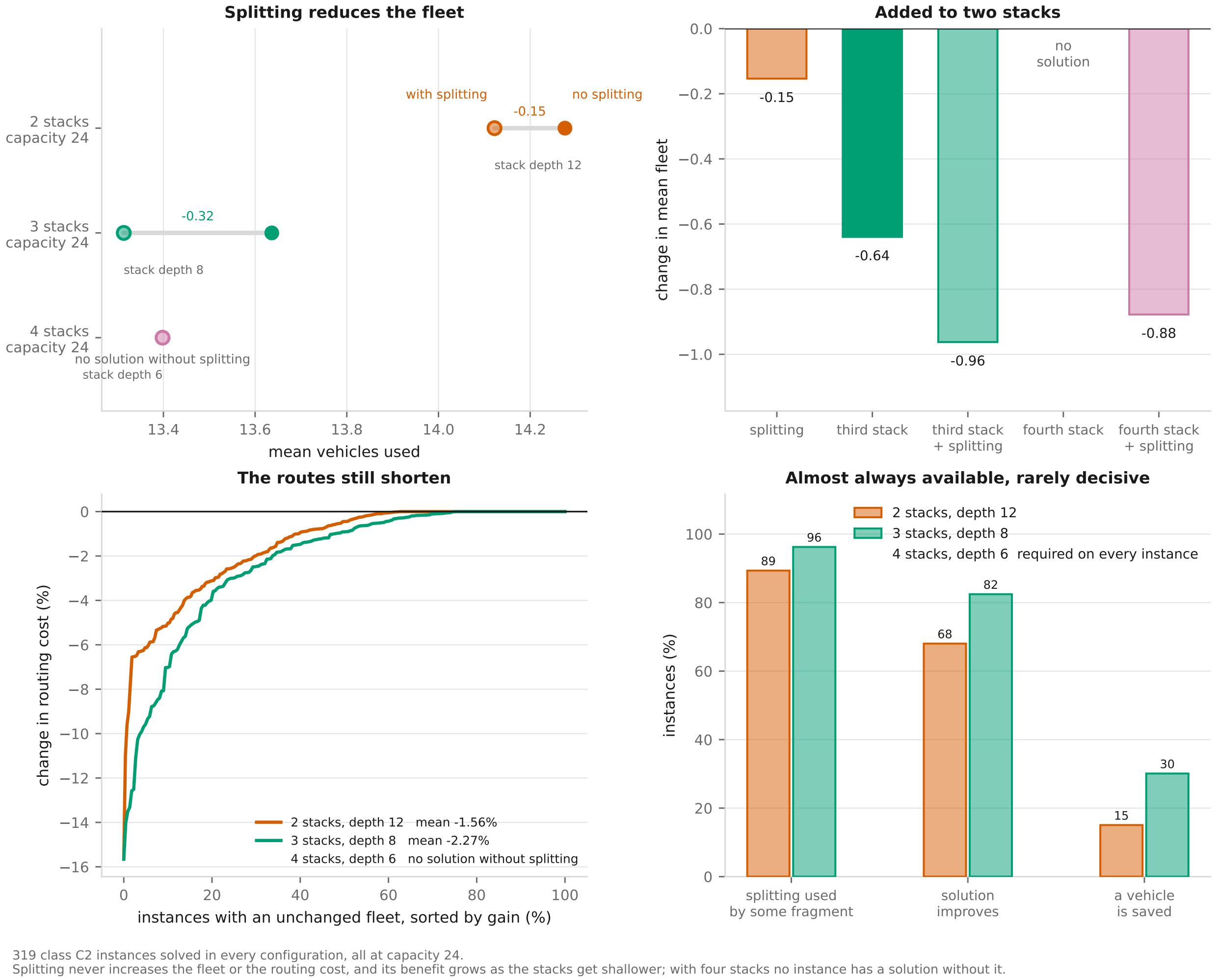

Figure 4: Two, three and four stacks, all at capacity 24. Top: mean fleet with and without splitting, and the change from adding each feature to two stacks without splitting. Bottom: routing cost where the fleet is unchanged, and how often splitting is available against how often it decides the answer. Four stacks have no solution without splitting, so the bottom panels have nothing to compare them with.

Fleet size is the headline objective and moves only in whole vehicles, so measuring splitting by it alone understates the effect. On the instances where the fleet does not move the routes still shorten, and the bottom left panel of Figure 4 sorts that change instance by instance. Splitting lowers routing cost by 1.56% over the 271 two-stack instances whose fleet is unchanged and by 2.27% over the 223 with three stacks, with a best case of 15.7% (Figure 4). The extra stack is worth more again on this measure, 4.12% at constant capacity, so splitting is the smaller of the two effects on routing cost as well as on fleet size.

Splitting is nonetheless available far more often than it is decisive, which the bottom right panel of Figure 4 separates out. A fragment reaches the splitting IP only after the ordinary PDPTWMS check has rejected it, so the fragments that IP accepts are exactly those that exist because splitting is allowed. At least one such fragment is generated on 89.3% of instances with two stacks and on 96.2% with three, yet they account for only 0.30% and 1.29% of all fragments generated. The solution improves on 68.0% and 82.4% of instances, and a vehicle is saved on 15.0% and 30.1%. Splitting is thus a small perturbation to the fragment set that is almost always available, often worth a little routing cost, and occasionally worth an entire vehicle.

Dividing the same capacity of 24 into four stacks produces the case the benchmark lacks. The depth falls to six, and every instance contains requests of size seven or eight that exceed it. Table 10 sets the four-stack solutions against those for two and three stacks at the same capacity, all three with splitting, over all 319 instances.

Table 10: Two, three and four stacks at a capacity of 24, all with splitting, over the 319 class C2 instances, all of which every configuration solves. Fleet and cost are compared with one stack fewer, cost only on the instances whose fleet is unchanged

| | *2 stacks* | *3 stacks* | *4 stacks* |
|---|---|---|---|
| Stack depth | 12 | 8 | 6 |
| Mean fleet | 14.12 | 13.31 | 13.40 |
| One vehicle fewer | – | 69.9% | 4.4% |
| One vehicle more | – | 0.0% | 12.9% |
| Cost at equal fleet | – | −4.98% | +0.58% |
| Non-dominated fragments (millions) | 7.4 | 16.3 | 17.2 |
| Stack assignment IPs solved | 64,611 | 350,471 | 1,771,890 |
| Share of time in fragment generation | 6.4% | 6.3% | 33.5% |
| Median solution time (s) | 0.14 | 0.20 | 0.45 |

Splitting makes the four-stack vehicle feasible, and a good deal better than a two-stack one. Against two stacks it removes a vehicle on 62.7% of instances and adds one on none. Against three stacks, though, the gain stops. The fourth stack adds a vehicle on 12.9% of instances and removes one on 4.4%, raising the mean fleet from 13.31 to 13.40, and where the fleet is unchanged it lengthens the routes by 0.58% on average. Adding stacks at a fixed capacity therefore buys accessibility until the stacks become too shallow to hold a request, and the fourth stack is past that point. Dividing 24 units four ways leaves room in a stack for two whole requests only when both are of size three, and spreads every request of size seven or eight over at least two stacks, so that its delivery constrains the order in every stack it occupies. Splitting lets the vehicle carry such requests at all, but it cannot give back the depth the fourth stack takes away. The results bear out the example of Section 3, in which three stacks carry a load that four of the same total capacity cannot, even with splitting.

The fourth stack also moves the computational effort into the loading problem. The fragment set grows only a little against three stacks, from 16.3 to 17.2 million, but the stack assignment IP is solved 1,771,890 times, five times as often as with three stacks and 27 times as often as with two, and fragment generation grows from 6.3% of the solution time to 33.5%. Step 2 of Table 7, which lets no two requests share a stack, becomes the first to resolve 40.0% of the calls, against 37.4% with three stacks and 13.2% with two, and the stack assignment IP resolves 7.1% of them against 1.5% with three stacks. The four-stack run is slower than the three-stack run on 274 of the 319 instances, by a factor of 1.79 in geometric mean, and its median instance takes 0.45 s against 0.20 s. It nonetheless solves all 319 instances within the time limit.

### *7.4. Splitting against compartments without LIFO*

Request splitting is not the only way to loosen the multi-stack policy. The PDPTWMS-H1 of Section 3 loosens it in a different direction. With zero handling time it keeps every request whole on one stack but drops the LIFO order entirely, so a fragment is accepted whenever the requests on board can be assigned to stacks within capacity. The two are alternative relaxations of the same problem, neither contained in the other, and the same method solves both, so what each is worth can be measured on the same instances against the same PDPTWMS baseline.

Table 11 does that over the a280 class C2 instances, at a vehicle capacity of 24 with two and with three stacks. Its left half measures each relaxation against the PDPTWMS it relaxes; its right half against the PDPTW of the paragraphs below. The a280 family is used because it is the hardest of the five and the one Section 6 reports in detail. Both relaxations are compared against the PDPTWMS solution of the same instance in the same configuration, over the 98 instances solved in every configuration, including the PDPTW of the paragraphs below.

Table 11: What each relaxation is worth, against the PDPTWMS it relaxes and against the PDPTW it reaches toward, over the 98 a280 class C2 instances solved in every configuration

| | | Against the PDPTWMS | | | Against the PDPTW | |
|---|---|---|---|---|---|---|
| *Configuration* | *Loading policy* | *Fleet* | *Saves a vehicle* | *Cost at equal fleet* | *Fleet above* | *Reaches it* |
| 2 stacks, depth 12 | PDPTWMS | – | – | – | +1.163 | 17.3% |
| | No LIFO (H1) | −0.969 | 77.6% | −9.84% | +0.194 | 80.6% |
| | Splitting | −0.112 | 11.2% | −2.61% | +1.051 | 19.4% |
| 3 stacks, depth 8 | PDPTWMS | – | – | – | +0.888 | 30.6% |
| | No LIFO (H1) | −0.286 | 28.6% | −3.63% | +0.602 | 43.9% |
| | Splitting | −0.541 | 48.0% | −5.17% | +0.347 | 65.3% |

The ordering reverses with the number of stacks, and the reversal is the point. With two stacks the LIFO restriction is worth far more than the indivisibility of a request. Dropping it removes a vehicle on 77.6% of the 98 instances and two or more on 17, against 11.2% and none for splitting. Instance by instance the H1 solution is at least as good as the split one everywhere, and strictly better on 96 of the 98. With three stacks splitting overtakes it, saving a vehicle on 48.0% of instances against 28.6%, and is strictly better on 80 of the 98.

Stack depth explains the exchange. Two stacks of depth 12 are deep enough that little is gained by dividing a request between them, but deep stacks bury many items, so LIFO is what binds and releasing it is worth a great deal. Three stacks of depth 8 bury less, which costs H1 most of its advantage, while the shallower stack is exactly what makes a request awkward to place, which is where splitting earns its gain. The two relaxations are therefore not competing accounts of the same story. They address different constraints, and which one matters depends on how the same capacity is divided.

Both relaxations can also be measured against the point they are reaching toward. Removing the stack structure altogether leaves the PDPTW, in which a vehicle is a single hold of the same total capacity and only that capacity binds. Our method solves it unchanged, since a fragment whose requests fit in the vehicle is then always loadable, and the optimum it returns is a lower bound on all three of the problems above. The right half of Table 11 measures each configuration against that floor, on the same 98 instances, on which the PDPTW needs a mean fleet of 10.31 vehicles. Where a configuration matches the floor's fleet its routes are still longer, by 3.86% for the three-stack PDPTWRSMS and 3.16% for two stacks without LIFO.

Read as a fraction of the distance the PDPTWMS has to travel, the reversal is sharper still. With two stacks, dropping LIFO closes 83.3% of the gap to the floor and splitting 9.6%; with three, splitting closes 60.9% and dropping LIFO 32.2%. The three-stack PDPTWRSMS is the configuration that comes closest to the floor of any we solve, needing 0.347 vehicles more than a vehicle with no stack structure at all and matching it exactly on 65.3% of instances. Dividing the capacity into three shallow stacks and allowing a request to be split across them therefore recovers most of what the multi-stack policy costs, without giving up the LIFO discipline that makes the policy worth having.

The comparison also bears on the cost of each relaxation in practice. Dropping LIFO is the larger change to the problem, and it is the harder one to justify operationally, since it assumes items can be dug out and replaced at no cost. Splitting assumes only that a request is made of pallets that need not travel together, which is the situation the PDPTWRSMS is proposed for, and at three stacks it delivers the larger saving of the two.

### *7.5. Route composition*

Recall from Section 5.1 that a route is a chain of fragments, joined wherever the vehicle becomes empty. This last question is one about the method rather than about the loading policy, so it is asked of the whole benchmark, over the 30,683 routes of the 1,910 instances solved to optimality within the hour. Table 12 reports how long those chains actually are.

Table 12: Number of fragments per vehicle route in the optimal solutions

| *Fragments in route* | *1 stack* | *2 stacks* | *3 stacks* | *All* |
|---|---|---|---|---|
| 1 | 69.91 | 87.18 | 87.94 | 79.72 |
| 2 | 28.54 | 12.38 | 11.46 | 19.30 |
| 3 | 1.53 | 0.42 | 0.61 | 0.96 |
| 4 | 0.01 | 0.02 | – | 0.01 |

Almost 80% of routes consist of a single fragment and a further 19% of two, so only 0.97% of routes (299 of 30,683) are built from three or more fragments. The proportion is highest with a single stack, where a vehicle empties more often and routes are therefore split more finely, and drops sharply once a second stack is available.

This has a direct consequence for the chaining constraints. The separation of Section 5.1 must in principle examine a whole chain, because travel time accumulates along it and three or more fragments can join with every adjacent pair feasible while the route as a whole is not. In practice such chains are rare, and a separation routine that tests only consecutive pairs is exact on 99.03% of the routes appearing in these optimal solutions, and over the 1,914 runs the lazy constraints were separated on 38 instances, for 58 cuts in all. The general procedure is nonetheless doing work as well as being necessary for correctness, since three of those cuts removed a chain of more than two fragments and one closed a cycle.

## 8. Conclusions

In this paper, we introduce a new practical variant of the PDPTWMS which we denote the PDPTWRSMS. In this variant we allow the size of a request to exceed the stack capacity and be only limited to the vehicle capacity, and we allow a request to be split among all stacks. This flexibility is very important for the transportation of pallets or other standardised containers and is essential when the size of any single request exceeds the size of a stack. Computational experiments show that applying this variant to the problem results in a significant cost reduction, especially as the number of stacks increases.

We present a new IP formulation to model the assignment of requests to stacks when the requests can be split across multiple stacks. We propose an algorithm to produce loading plans for the PDPTWRSMS which takes advantage of the problem characteristics to minimise the number of stack assignment IPs solved.

Splitting is not the only way to loosen the multi-stack policy, and Section 7.4 compares it with the alternative of keeping requests whole and dropping the LIFO discipline instead. Which is worth more depends on how the capacity is divided. Measured against the PDPTW, in which the vehicle has no stack structure at all, splitting recovers 9.6% of the available fleet reduction with two deep stacks and 60.9% with three shallow ones, while dropping LIFO recovers 83.3% and 32.2%. With three stacks and splitting the fleet is within 0.35 vehicles of a vehicle with no stacks whatever. Where the stacks are shallow, splitting therefore recovers most of what the multi-stack policy costs, and it does so without giving up the LIFO discipline that makes the policy worth having.

We modify the algorithm of Alyasiry et al. (2019) to solve the PDPTWMS. Computational results show that our method outperforms the branch-price-and-cut algorithm of Cherkesly and Gschwind (2022), the current state of the art, solving 1,910 of the 1,914 benchmark instances within an hour against 1,731, closing at least 179 that were previously unsolved, and the remaining four within eight hours. Mean solution times fall by factors of 15 to 155 on the two- and three-stack instances, and by about 250 on the single-stack instances of class C2 and 900 on those of class C1.

There is considerable scope for future work on the method of fragments and its application to similar problems. For large instances, the time to generate fragments is a large fraction of the solution time, and our implementation generates them on a single thread; the enumeration is naturally parallel and would benefit directly from being distributed across cores.

Future work could also consider extending the method of fragments to include additional loading constraints, such as incompatible requests or trucks with segmented storage (e.g. normal and refrigerated). Aerts-Veenstra et al. (2024) and Aerts-Veenstra et al. (2025) study constraints of this kind for compartments that impose no loading order, and combining them with the multi-stack policy would be a natural next step. Additionally, work is required to extend the method of fragments to cases where the number of possible fragments grows extremely large. A benchmark in which the number of stacks and the total vehicle capacity varied independently would also help, since in the present instances a third stack arrives together with a rise in capacity from 24 to 27 for class C2, and the two effects cannot be separated from the published configurations alone. Instances that exercise the separation of Section 5.1 more heavily would be valuable for the same reason. Chains needing the general procedure are rare in the present set, where the lazy constraints were separated on 38 of the 1,914 runs for 58 cuts in all, only four of which went beyond a consecutive pair, so the benchmark tests that part of the method very lightly.

## Appendix A. Detailed results

Table A.1: PDPTWRSMS for 3 stacks vehicle with capacity 24

| *Instance* (a280-C2) | $F_N$ | *IPs* | *Gap* | $\text{Sec}_{\text{IP}}$ | $\text{Sec}_f$ | *Sec.* | $Z^*$ | $S_{\text{Veh}}$ | $S_{\text{Dist}}$ |
|---|---|---|---|---|---|---|---|---|---|
| 141-w15-500-1000 | 5,494 | 232 | 5.625 | 0.0 | 0.1 | 0.3 | 1,913,172.72 | 0 | −5.18 |
| 151-w15-500-1000 | 4,489 | 164 | 4.127 | 0.0 | 0.1 | 1.3 | 2,014,859.89 | −1 | −3.07 |
| 141-w15-1000-1200 | 13,485 | 103 | 1.870 | 0.0 | 0.1 | 0.9 | 1,712,684.46 | 0 | −2.90 |
| 151-w15-1000-1200 | 16,849 | 284 | 2.169 | 0.0 | 0.2 | 2.0 | 1,614,121.42 | −1 | 2.82 |
| 141-w15-1500-2000 | 148,838 | 1,653 | 7.636 | 0.5 | 2.0 | 60.8 | 1,211,520.72 | 0 | −13.51 |
| 151-w15-1500-2000 | 357,062 | 5,046 | 2.750 | 1.3 | 5.3 | 141.1 | 1,212,957.05 | −1 | 5.09 |
| 141-w30-500-1000 | 15,230 | 572 | 4.765 | 0.1 | 0.3 | 4.4 | 1,612,630.64 | −1 | −3.23 |
| 151-w30-500-1000 | 19,528 | 578 | 2.987 | 0.1 | 0.3 | 3.4 | 1,612,876.95 | −2 | −3.53 |
| 141-w30-1000-1200 | 53,409 | 446 | 6.547 | 0.1 | 0.5 | 5.1 | 1,411,489.02 | 0 | −4.22 |
| 151-w30-1000-1200 | 50,734 | 1,015 | 7.267 | 0.3 | 0.8 | 24.9 | 1,512,163.17 | 0 | −2.34 |
| 141-w30-1500-2000 | 821,773 | 13,868 | 8.272 | 5.7 | 17.3 | 567.3 | 1,210,948.39 | −1 | −2.78 |
| 151-w30-1500-2000 | 790,654 | 6,927 | 6.682 | 2.3 | 13.1 | 129.5 | 1,212,116.37 | −1 | −2.10 |
| 141-w45-500-1000 | 26,273 | 1,073 | 2.097 | 0.2 | 0.6 | 6.8 | 1,613,025.82 | −2 | 3.91 |
| 151-w45-500-1000 | 66,323 | 1,857 | 2.983 | 0.4 | 1.4 | 30.2 | 1,612,495.11 | −1 | −5.94 |
| 141-w45-1000-1200 | 87,663 | 2,164 | 7.725 | 0.6 | 1.9 | 34.7 | 1,311,115.67 | −1 | −3.54 |
| 151-w45-1000-1200 | 199,346 | 3,689 | 3.232 | 1.1 | 4.2 | 55.3 | 1,312,362.31 | −1 | 2.44 |
| 141-w45-1500-2000 | 1,153,477 | 29,589 | 4.494 | 10.2 | 33.6 | 413.5 | 1,211,375.72 | −1 | −5.21 |
| 151-w45-1500-2000 | 3,617,549 | 92,854 | 7.591 | 46.4 | 115.5 | 3,086.7 | 1,211,173.69 | −1 | −1.76 |